\documentclass[12pt, a4paper]{article}
\usepackage{geometry}

\usepackage{comment}
\usepackage{enumerate}
\usepackage[T1]{fontenc}
\usepackage{latexsym}
\usepackage{amssymb}
\usepackage{amsmath}
\usepackage{amsthm}
\usepackage{dsfont}
\usepackage{mathrsfs}
\usepackage{blindtext}
\usepackage{ucs}
\usepackage{subfigure}
\usepackage{multirow}
\usepackage{graphicx}
\usepackage{color}
\usepackage{xcolor}
\usepackage[pdftex,bookmarks=blue,citecolor=blue,urlcolor=cyan,colorlinks,linkcolor=blue]{hyperref}            

\newtheorem{thm}{{\sc Theorem}}[section]
\newtheorem{prop}[thm]{{\sc Proposition}}
\newtheorem{lem}[thm]{{\sc Lemma}}

\theoremstyle{definition} 
\newtheorem{defn}[thm]{Definition} 

\theoremstyle{Example} 
\newtheorem{ex}[thm]{Example}  
\newcounter{saveeqn}

 \numberwithin{equation}{section}

\title{Affine AP-Frames and Stationary Random Processes}
\date{}
\begin{document}
\author{{\bf {\footnotesize BY}}\\ 
{\bf {\normalsize Hern\'{a}n D. Centeno}}\\
{\footnotesize AND} {\bf {\normalsize Juan M. Medina}}}
\maketitle
\noindent
{\it Abstract.} {\small 
It is known that, in general, an affine or Gabor AP-frame is an $L^2(\mathbb{R})$-frame and conversely. 
In part as a consequence of the Ergodic Theorem, we prove a necessary and sufficient condition for an affine (wavelet) system $\mathcal{A}=\{a^{j/2} \psi_{j,k}(t):=a^{-j/2} \psi (a^{-j} t -k) :j\in\mathbb{Z}, k\in\mathbb{K}:=b\mathbb{Z}\}$ to be an affine AP-Frame in terms of Gaussian stationary random processes expanding in this way what we have done recently for Gabor systems.    
Likewise, we study a connection between the decay of the associated stationary sequences $\{\langle{X,\psi_{j,k}}\rangle : k\in\mathbb{K}\}$ for each $j\in\mathbb{Z}$, and a smoothness condition on a Gaussian stationary random process $X=(X(t))_{t\in\mathbb{R}}$.
}\\[0.2cm]
\\[0.2cm]
\maketitle

\setcounter{section}{0} 
\section{Introduction}

The most important spaces of almost periodic functions are non-separable and therefore they cannot admit countable frames. 
However, several authors \cite{galindo, Kim2009, Partin} introduced the related concept of \emph{AP-frame} for \emph{Gabor} and \emph{affine} systems and proved conditions under which frame-type inequalities are still possible as estimators of the norm in these spaces.
Under rather mild conditions it results that, for these systems, to be an AP-frame and a $L^2(\mathbb{R})$-frame are equivalent notions. 
In a previous work \cite{CenMed22} the authors of this paper have seen how some of these results for Gabor systems $\mathcal{G}$ are connected or equivalent to other involving stationary random processes. This connection arises naturally, since stationary random processes with discrete spectrum, i.e. with a discrete spectral measure, are closely related to almost periodic functions. In contrast, the same problem when considering stationary random processes with an arbitrary spectral measure, deserves a detailed study. 
In this general case, this relationship, in part, relies in some classic results from Ergodic theory (see Subsection \ref{StoProc}).  
In Section \ref{AP-frames-processes} are obtained several results for affine systems $\mathcal{A}$ analogous to the ones for Gabor systems $\mathcal{G}$ in \cite{CenMed22}. 
By the different nature of these systems (for instance a Fourier transformed Gabor system $\widehat{\mathcal{G}}$ is again a Gabor system but generally this does not happen for affine systems) some results obtained here have required alternative techniques and approaches, some of them adapted or inspired from those in \cite{Kim2009}. 
Apart of giving an alternative characterization of affine systems which constitute a frame, we will see that the same techniques are useful for retrieving some satistical information about the smoothness of the involved random processes. We note that, despite Wavelet analysis is a common tool in the study of smoothness in terms of function spaces, there is no known reference of it in the aforementioned articles in the context of  \emph{AP-frames}. This notion will be introduced here, apparently, for the firs time.
In fact, for a Gaussian stationary random process $(X(t))_{t\in\mathbb{R}}$ its associated $\alpha$-th derivative process $(D^\alpha X (t))_{t\in\mathbb{R}}$, with $\alpha>0$, is meaningful if it can be seen as an admissible linear time invariant filtering of $X$ (Subsection \ref{StoProc} equation \eqref{filterfreq}) with a filter in $L^2(\mathbb{R},\mu)$, where $\mu$ is the spectral measure associated to $X$. 
In this case $D^\alpha X$ is also a Gaussian stationary random process.
In Section \ref{fractional} it is shown that this filter representation, or equivalently the existence of the $\alpha$-th derivative, is equivalent to a pair of (essentially different) conditions. 
First, to the decay of the affine coefficients $\langle X,\psi_{j,k}\rangle$ associated to X ($0<\alpha< 1$), and second, to some sort of Sobolev-type integrability condition involving: a singular integral of the covariance function of $X$ ($0<\alpha< 1$) as well as a singular integral of a centered difference of the original process $X$ ($0<\alpha< 2$). 
The frame property on both, the affine system $\mathcal{A}$ and the system $I^\alpha\mathcal{A}$ (the system $\mathcal{A}$ transformed by the Riesz potential $I^\alpha$) plays an important role in the first equivalence while in the second one are involved techniques of Gaussian stationary random processes together with those related to singular integrals.

\section{Preliminaries}

We shall thoroughly recall the terminology and results already introduced,, respectively, in the preliminary results sections of the articles \cite{CenMed22, Med}. 
In order to ease the reading of the present article we make a brief summary of these also skipping several details.
 
\subsection{Some background on function spaces and the Fourier transform}\label{func}

As usual, for $p\in [1,\infty) $, we will denote the classical Lebesgue function spaces with $L^p(\mathbb{R})$. 
When $p=2$, we consider $L^2(\mathbb{R})$ endowed with its usual inner product. 
With some abuse, we shall use the same notation when this integral is well defined for functions which not necessarily belong to $L^2 (\mathbb{R})$.
The \emph{Fourier Transform} of $f\in L^1(\mathbb{R})$ is given by:
$$\mathcal{F}f(\lambda)=\widehat{f}(\lambda)=\int_{\mathbb{R}} f(x) e^{-i\lambda x}\,dx\,.$$
Analogously, if $\widehat{f}$ is integrable, $f$ can be recovered by the inverse Fourier Transform, 
$$(\widehat{f})^{\vee}(x)=\frac{1}{2\pi} \int_{\mathbb{R}}\widehat{f}(\lambda) e^{i\lambda x}\,d\lambda\,.$$
By a density argument the Fourier Transform can be defined for $f\in L^2(\mathbb{R})$. 
In fact, in this case, one has the \emph{Plancherel identity:}
$$\Vert\widehat{f}\,\Vert^2 _{L^2(\mathbb{R})}= 2\pi \Vert f\Vert^2 _{L^2(\mathbb{R})}\,$$
expressing the fact that the Fourier Transform, over $L^2 (\mathbb{R})$, is a unitary map. 
Fourier transforms are defined for other classes of measures or functions. 
For more details see for example \cite{Katz,Pinsky}. 
If $\mathcal{G}$ denotes any family or subset of functions for which its Fourier transform is defined, we will denote its image as $\widehat{\mathcal{G}}=\mathcal{F}(\mathcal{G})=\{\widehat{f}:\,f\in\mathcal{G}\}$. 
We shall introduce another class of functions. 
Recall that $AP(\mathbb{R})$ coincides with the uniform norm closure of the space of trigonometric polynomials $\sum_jC(j)e^{i\lambda_j t}$ with $\lambda_j\in \mathbb{R}$ and $C(j)\in\mathbb{C}$. 
In $AP(\mathbb{R})$ one can introduce the inner product:
\begin{equation}\label{pi}
\langle f,g\rangle_{AP(\mathbb{R})}=
\lim\limits_{T\longrightarrow\infty}\frac{1}{2T} \int_{-T}^T f(t)\overline{g(t)}\,dt\,.
\end{equation}
The norm induced by this inner product makes $AP(\mathbb{R})$ a non-complete and non-separable inner-product space. 
Now, we recall the Hilbert space $B^2 (\mathbb{R})$ of \emph{Besicovitch almost periodic functions} containing $AP(\mathbb{R})$. 
Let $f\in L^2 _{loc}(\mathbb{R})$ we can define the semi-norm:
$$\left\|f\right\|_{b,\,2} = \left({\limsup\limits_{T\longrightarrow\infty}
\frac{1}{2T} \int_{-T}^T |f(t)|^2 d t}\right)^{\frac{1}{2}}\,.$$
A function $f\in\,L^2_{loc}(\mathbb{R})$ is called \emph{Besicovitch} almost periodic, if for every $\varepsilon>0$ there exists $g\in AP(\mathbb{R})$ such that $\left\|f-g\right\|_{b,\,2} <\varepsilon$. 
It is possible to turn $B^2 (\mathbb{R})$ into a Hilbert space. 
First one can define an equivalence relation on $B^2 (\mathbb{R})$ in the following way: $f\equiv g$ if and only if $\left\|f-g\right\|_{b,\,2}=0$. 
The norm of $[f]\in B^2 (\mathbb{R})/\equiv $ is given by $\left\|[f]\right\|_{B^2(\mathbb{R})}:= \left\|f\right\|_{b,2}$. 
In particular if $f\in L^2 _{loc} (\mathbb{R})$ is such that $\left\|f\right\|_{b,\,2}=0$ then $f\equiv 0$. 
Finally, one can prove that $B^2 (\mathbb{R})/\equiv$ is complete. See e.g. \cite{Cor}, 
or for a detailed exposition of this construction in LCA groups see e.g. \cite{Lenz2020} p.39. 
With some abuse, if there is no confusion, we will write $f$ instead of $[f]$ the equivalence class of $f$ and $B^2 (\mathbb{R})$ for $B^2 (\mathbb{R})/\equiv$. 
Moreover, it can be proved that the inner product $\langle{\,.\,,\,.\,}\rangle_{B^2 (\mathbb{R})}$ in $B^2 (\mathbb{R})$ coincides with equation \eqref{pi} for any $f,g \in  B^2 (\mathbb{R})$. 

In $B^2 (\mathbb{R})$, the complex exponentials $(e^{i\lambda t})_{\lambda\in\mathbb{R}}$ form a complete orthonormal basis and the following analogue, due to Wiener \cite{Katz}, of Plancherel identity holds:
\begin{equation}\label{wien1}
\Vert f\Vert_{B^2(\mathbb{R})} =\Vert C(f)\Vert_{L^2(\mathbb{R}, dc)}\,,
\end{equation}
where $C(f)(\lambda)=\lim\limits_{T\longrightarrow\infty} \frac{1}{2T}\int_{-T}^T f(t) e^{-i\lambda t}\,dt$ denotes the \emph{Bohr Transform} of $f$ and $c$ denotes the counting measure. 
If $b>0$ and $\mathbb{K}=b\mathbb{Z}$, in an analogous way, one can introduce the space of almost periodic sequences $AP(\mathbb{K})$ (or over any subgroup $\alpha \mathbb{Z}$). 
Indeed an almost periodic sequence $(x(k))_{k\in\mathbb{K}}$ can be characterized as the restriction over $\mathbb{K}$ of functions in $AP(\mathbb{R})$. 
This space can be endowed with the inner product  
$$\langle x,y\rangle_{AP(\mathbb{K})}=\lim\limits_{N\longrightarrow\infty}\frac{1}{2N+1} \sum\limits_{k\in\mathbb{K}(N)} x(k)\overline{y(k)} \,,$$
where $\mathbb{K}(N)=\{n b\in\mathbb{K}:-N\leq n\leq N\,\}$. 
Again, these spaces admit a completion $B^2(\mathbb{K})$ and $(e^{i\lambda k})_{\lambda\in\mathbb{T}}$ is a non-countable orthonormal basis of it. 
Here, Wiener's formula takes the form
\begin{equation}\label{wien2}
\Vert x\Vert_{B^2(\mathbb{K})} =\Vert C(x)\Vert_{L^2(\mathbb{T}, dc)}\,,
\end{equation}
with $C(x)(\lambda)=\lim\limits_{N\longrightarrow\infty}\frac{1}{2N+1} \sum\limits_{k\in\,\mathbb{K}(N)} x(k) e^{-i\lambda k}$. 

In a similar way to the definition of $\mathbb{K}$ and $\mathbb{K}(N)$, we introduce $\mathbb{D}=\frac{2\pi}{b} \mathbb{Z}$. 
The respective finite subsets $\mathbb{D}(N)$ and $\mathbb{Z}(N)$ are defined analogously to $\mathbb{K}(N)$. 
To avoid repetitions, we summarize some remarkable facts about $B^2 (\mathbb{R})$, analogous results hold for $\mathbb{K}$:
\begin{enumerate}[1.]
\item (Riesz-Fischer property) Let $C:\mathbb{R}\longrightarrow \mathbb{C}$. 
Then there exists a unique $f \in\,B^2 (\mathbb{R}) $ (as an equivalence class) such that $C(\lambda)=\langle{f,\, e^{i.\lambda(\,.\,)}\,}\rangle_{B^2 (\mathbb{R})}$ if and only if $\left\|\,C\,\right\|_{L^2({\mathbb{R}}, dc)}<\infty$. 
In this case $f=\sum\limits_{\lambda\in \mathbb{R}} C(\lambda) e^{i.\lambda(\,.\,)}\,$, where the convergence is in the $B^2 (\mathbb{R})$-norm.
\item Let $f\in AP(\mathbb{R})$, $\tau \in \mathbb{R}$ and let $T_\tau f= f(\,.\,+\tau) \in AP(\mathbb{R})$ be the translation of $f$ by $\tau$. 
Then $\left\|\,f\,\right\|_{B^2(\mathbb{R})}= \left\|\,T_\tau f\,\right\|_{B^2(\mathbb{R})}$ and hence $T_\tau$ extends uniquely to an isometry on $B^2 (\mathbb{R})$. 
Here, with some abuse we will shall also denote this extension by $T_\tau$. 
Moreover, one can define the \emph{deterministic auto-correlation} of $f \in B^2 (\mathbb{R})$ at $\tau \in \mathbb{R}$ by:  
\begin{equation}\label{autocorr}
\rho_f (\tau):=\langle{\,f\,,\,T_\tau f\,}\rangle_{B^2 (\mathbb{R})}\,.
\nonumber
\end{equation}
If $f\in B^2 (\mathbb{R}) $ then $\rho_f \in AP(\mathbb{R})$. 
In fact, from \eqref{wien1} one can deduce that 
$$\rho_f(\tau)= \sum\limits_{\lambda\in\mathbb{R}}|C(\lambda)|^2 e^{-i.\lambda \tau}\,,$$
and therefore $\rho_f$ is the Fourier transform of the discrete measure $\nu=\sum\limits_{\lambda}|C(\lambda)|^2 \delta_{\lambda}$.
\end{enumerate}

\subsection{Hilbert frames and affine AP-frames}

\subsubsection{Hilbert spaces and frame sequences}

We are interested in certain frame sequences and some of their properties. 
Let $\mathcal{H}$ be a complex Hilbert space with inner product $\langle\,\cdot\,,\,\cdot\,\rangle_{\mathcal{H}}$ (linear in the first variable and anti-linear in the second one) and $\mathbb{J}$ a countable set of indexes. 
We recall the main definitions and for a comprehensive treatment of these topics we refer the reader to e.g. \cite{Chris}.  
A sequence $\{f_j\}_{j\in\mathbb{J}} \subset \mathcal{H}$ is a \emph{frame}, or \emph{frame sequence}, for $\mathcal{H}$ if there exist constants $0<A\leq B$ such that for every $f\in \mathcal{H}$:
\begin{equation}\label{frame}
A \left\|f\right\|^2 \leq \sum\limits_{j \in \mathbb{J}} |\langle{f,f_j}\rangle_{\mathcal{H}} |^2 \leq B \left\|f\right\|^2 \,.
\end{equation}
Sometimes, a frame is thought informally as a redundant, stable and complete system in $\mathcal{H}$. 
In other words, any vector belonging to $\mathcal{H}$ has a, not necessary unique, unconditionally convergent expansion with respect to $\{f_j\}_{j\in\mathbb{J}}$. \hfill\\
Finally, Let $T:\mathcal{H}\longrightarrow \mathcal{H}$ be a bounded linear operator. 
We say that $T$ is \emph{positive} if $\langle{Tx,x}\rangle_{\mathcal{H}}\geq 0$ for all $x\in\mathcal{H}$ and $T$ is \emph{positive definite} if $\langle{Tx,x}\rangle_{\mathcal{H}} > 0$ for all $x\in\mathcal{H}\smallsetminus \{0\}$. 
If $T$ is positive then it is self-adjoint \cite[p.168]{Kantor}.

\subsubsection{Affine systems and AP-frames}\label{TF}

Since $B^2(\mathbb{R})$ is a non-separable Hilbert space, this prompts an alternative definition to \eqref{frame}. 
We begin with some auxiliary statements. 

\begin{defn}\label{AS} 
Let $\psi:\mathbb{R}\longrightarrow\mathbb{C}$ be a Borel measurable function (the \textit{mother wavelet}) and $a>1$, $b>0$ fixed constants. 
We define 
\begin{equation}\label{AL1Norm}
\psi_{j,k}(t):=
a^{-j}\psi(a^{-j} t -k),\hspace{0.5cm} j\in\mathbb{Z},k\in\mathbb{K}:=b\mathbb{Z}\,,
\nonumber
\end{equation} 
and the affine (wavelet) system generated by time-scale shifts of $\psi$ as
\begin{equation}\label{ASdef}
\mathcal{A}=\mathcal{A}(\psi,a,b)=
\{a^{j/2}\psi_{j,k}(t):=
a^{-j/2}\psi(a^{-j} t -k),\,j\in\mathbb{Z},\ k\in\mathbb{K}\}\,.
\nonumber
\end{equation}
\end{defn}
A crucial topic in wavelet analysis is to characterize when $\mathcal{A}\subset L^2(\mathbb{R})$ is a frame for $L^2(\mathbb{R})$. 
Recalling that the Fourier Transform is a unitary map, $\mathcal{A}$ is a frame for $L^2(\mathbb{R})$ if and only if $\widehat{\mathcal{A}}$ is also a frame. 
Note that $\widehat{\mathcal{A}}$ may not be an affine system:
\begin{equation}\label{ASF}
\widehat{\mathcal{A}}=\{a^{j/2}\widehat{\psi}_{j,k}(\lambda)=a^{j/2}\,\widehat{\psi}(a^{j}\lambda)e^{-ika^{j}\lambda},\,j\in\mathbb{Z},\,k\in\mathbb{K}\}\,.
\end{equation}
Hereafter, $a\in\,\mathbb{N}$ and $b>0$. 
Nevertheless, for clarity in the exposition, without loss of generality the reader may assume that  $b=1$ (See e.g. \cite{RonShen97} p.419). 
Motivated in part by some interpolation problems for stationary random process we shall work preferable with (\ref{ASF}). 
At this point we introduce the affine product: 
\begin{equation}\label{eq-affprod}
[\lambda,\lambda']_\psi:=\sum_{j\geq\kappa(\lambda-\lambda')}\widehat{\psi}(a^j\lambda)\overline{\widehat{\psi}}(a^j\lambda')\hspace{0.5cm}\lambda,\lambda'\in\mathbb{R}\,,
\nonumber
\end{equation}
and then, at least formally, for each $\lambda\in\mathbb{R}$ the infinite matrices 
$({G}_j(\lambda)(d,d'))$, 
$(G(\lambda)(d,d'))$, 
$(\mathds{G}(\lambda)(q,q'))$ and 
$(\mathds{G}_j(\lambda)(q,q'))$, for $d,d'\in \mathbb{D}$ and $q,q'\in Q$, defined by 
\begin{align}
{G}_j(\lambda)(d,d')&=\widehat{\psi}(a^j(\lambda+d))\overline{\widehat{\psi}}(a^j(\lambda+d'))\mathbf{1}_{\mathbb{D}}(a^j(d-d'))\hspace{0.5cm}j\in\mathbb{Z},\nonumber \\ 
G(\lambda)(d,d')&=[\lambda+d,\lambda+d']_\psi,\label{Gram}  \\
\mathds{G}_j(\lambda)(q,q')&=\widehat{\psi}(a^j(\lambda+q))\overline{\widehat{\psi}}(a^j(\lambda+q')\mathbf{1}_{\mathbb{D}}(a^j(q-q'))\hspace{0.5cm}j\in\mathbb{Z},\nonumber \\
\mathds{G}(\lambda)(q,q')&=[\lambda+q,\lambda+q']_\psi,\hspace{0.5cm}q,q'\in Q:=\bigcup_{j\in\mathbb{Z}}{a^{-j}}\mathbb{D},\label{Grammodif}
\end{align}
where $\kappa :\mathbb{R}\longrightarrow\mathbb{Z}$ is the $\frac{2\pi}{b}$-$a$-\textit{adic valuation function} given by
\begin{equation}\label{valuation}
\kappa (\lambda)=\inf\{j\in\mathbb{Z}:a^j\lambda \in \mathbb{D}\}.
\nonumber
\end{equation}
Thus, in first place, $\kappa (0)=-\infty$, and $\kappa(\lambda)=+\infty$ unless $\lambda$ is a $\frac{2\pi}{b}$-$a$-adic integer:
\begin{equation}
\kappa(\lambda)=m\hspace{0.5cm}\Leftrightarrow\hspace{0.5cm}a^m\lambda =d\hspace{0.5cm}\Leftrightarrow\hspace{0.5cm}\lambda=a^{-m}d\hspace{0.5cm}\textnormal{for some }d\in\mathbb{D}\,.
\nonumber
\end{equation} 
Second, $\mathbf{1}_{\mathbb{D}}(a^j(q-q'))=1$ if and only if $a^j(q-q')\in\mathbb{D}$ if and only if  $q-q'\in a^{-j}\mathbb{D}$. 
When $a\notin\mathbb{Q}$ we also have $q-q'\in a^{-j}\mathbb{D}$ if and only if $q,q'\in a^{-j}\mathbb{D}$. 

Instead of \eqref{Gram}, sometimes is preferable to work with \eqref{Grammodif}.
As pointed out in \cite{Kim2009}, for each $j\in\mathbb{N}$ the matrix $G(a^j\lambda)$ is obtained from $\mathds{G}(\lambda)$ by deleting all columns and rows indexed by $Q\smallsetminus a^{-j}\mathbb{D}$. 
This eliminates the redundancy which arises in $G(\lambda)$ being a submatrix of $G(a\lambda)$.

If, for instance, $\mathcal{A}$ is an $L^2(\mathbb{R})$-frame \cite{RonShen97}, then a.e. $\lambda\in\mathbb{R}$ we have
\begin{equation}\label{sumagrammianitos}
G(\lambda)(d,d')=\sum\limits_{j\in\mathbb{Z}}G_j(\lambda)(d,d')
\hspace{0.4cm}\textnormal{and}\hspace{0.4cm}
\mathds{G}(\lambda)(q,q')=\sum\limits_{j\in\mathbb{Z}}\mathds{G}_j(\lambda)(q,q')\,.
\end{equation}
Now we focus on $\mathds{G}$ but the same holds true for $G$ (in fact most of the results were originally proved for $G$). 
It is possible to characterize affine systems which are frames of $L^2(\mathbb{R})$ in terms of $\mathds{G}(\lambda)$ acting as a linear operator $\mathds{G}(\lambda): \ell^2(Q) \longrightarrow \ell^2 (Q)$.
In fact, a necessary condition is: for each $q\in Q$, $\mathds{G}(\cdot)(0,q)\in C(\mathbb{R}\setminus \{0,q\})$, and moreover, the following holds, cf. \cite[Corollary 3.2]{Kim2009}.

\begin{thm}\label{RS}
Let $\mathcal{A}\subset L^1(\mathbb{R})$ be an affine system and the associated matrix $\mathds{G}(\lambda)$ be defined as in (\ref{Grammodif}). 
Then: $\mathcal{A}$ (equivalently $\widehat{\mathcal{A}}$) is an $L^2(\mathbb{R})$-frame with constants $0<A\leq B$ if and only if:
\begin{equation}\label{RSpos}
A\left\|x\right\|^2 _{\ell^2(Q)}\leq\langle{\mathds{G}(\lambda)x,x}\rangle_{\ell^2(Q)}\leq B\left\|x\right\|^2 _{\ell^2(Q)}\,,
\end{equation}
for almost all (with respect the Lebesgue measure) $\lambda\in\mathbb{R}$ and all $x\in\,\ell^2(Q)$.
\end{thm}
For more details on linear operators defined by infinite matrices see \cite{Kantor}.    
As $B^2(\mathbb{R})$ is non-separable, therefore no countable collection of elements of its dual space can be total on it. 
However, by a suitable averaging process one can estimate the $B^2(\mathbb{R})$-norm. 
In fact, this motivates the following definition for affine systems \cite{ galindo,Kim2009,Kim2013,Partin}.

\begin{defn} 
Let $\mathcal{A}\subset L^1(\mathbb{R})$ be an affine system. 
We say that $\mathcal{A}$ is an affine AP-frame if there exist constants $0<A\leq B$ such that for every $f\in B^2(\mathbb{R})$ that satisfies $C(f)(0)=0$, it holds the following inequalities:
\begin{equation}
A\left\|f\right\|^2_{B^2(\mathbb{R})}\leq\sum\limits_{j\in\mathbb{Z}}\left\|\langle{f,\psi_{j,k}}\rangle\right\|^2_{B^2(\mathbb{K})}\leq B\left\|f\right\|^2_{B^2(\mathbb{R})}\,.
\nonumber
\end{equation}
\end{defn}
In \cite{Kim2009}, with additional mild conditions (as in Theorem \ref{RS} above), it is proved that an affine system is an $L^2(\mathbb{R})$-frame if and only if is an affine $AP$-frame.

\subsection{Stationary processes}\label{StoProc}
 
This brief background follows closely \cite{Pour,Ro,Roz}. 
Let $(\Omega,\mathcal{F},\mathbf{P})$ be a probability space and let $X=(X(t))_{t\in\mathbb{R}}\subset L^2(\Omega,\mathcal{F},\mathbf{P})$ be a complex, mean square continuous wide sense stationary (\emph{w.s.s.} for short) random process, i.e. $X$ verifies the following three conditions, for all $t,s\in\mathbb{R}$:
\begin{equation}\label{estat}
\mathbf{E}(X(t))=0,\hspace{0.3cm}
\mathbf{E}(X(t)\overline{X(s)})=R_X(t-s)\hspace{0.3cm}\textnormal{and}\hspace{0.3cm}
\lim_{t\rightarrow 0}\mathbf{E}|X(t)-X(0)|^2=0\,.
\end{equation}
Where $\mathbf{E}(\,.\,)$ denotes the expected value operator. 
In general, if $X$ is w.s.s. complex random process, we assume that $X(t)=X_1 (t) + i X_2 (t)$ for all $t\in \mathbb{R}$. 
Where $X_i$, $i=1,2$ are two stationary (cross)correlated \emph{real} w.s.s. stationary random processes. 
If $X$ is \emph{Gaussian} and \emph{complex} we shall impose, in addition to \eqref{estat}, the condition:
\begin{equation}\label{complwss}
\;\; \mathbf{E}(X(t)X(s))=0\;\;\textrm{for all}\;\; t,s\;\in \mathbb{R}\,.
\end{equation}
Gaussian complex random processes or vectors verifying condition \eqref{complwss} are said to be \emph{circular}. 
Mean square continuity, in addition, implies the existence of an equivalent random process $X^*$ (i.e. a process such that $\mathbf{P}(X(t)\neq X^*(t))=0$ for all $t$) which is measurable with respect to the completion of the product $\sigma$-algebra $\mathcal{F}\otimes\mathcal{B}(\mathbb{R})$. 
From this, we will usually assume, with no further mention of this fact, that we are working with such equivalent process when considering some operations on $X$ as function of $t$. 
A stronger notion is (strict) stationarity, i.e. if the shifted families $X^T=(X(t+T))_{t\in\mathbb{R}}$ have the same finite distributions as $X$ for all $T\in\mathbb{R}$. 
A strictly stationary process is w.s.s. but the converse is not always true. 
If $X$ is \emph{Gaussian} both notions coincides.    
If $X$ is a w.s.s. random process it is known by Bochner's Theorem that there exists a finite symmetric Borel measure $\mu$, \emph{the spectral measure}, such that the covariance function can be written in the following way:
$$R_X(t-u)=\mathbf{E}(X(t)\overline{X(u)})=\int_{\mathbb{R}}e^{i\lambda(t-u)}\,d\mu (\lambda)\,,$$
Conversely, if $\mu$ is a finite Borel measure, there exists a w.s.s random process with $\mu$ as its spectral measure. 
The \emph{spectrum of} $X$ is the support of $\mu$. 
In the case that $\mu$ is absolutely continuous with respect to the Lebesgue measure, then there exists its \emph{Radon-Nykodim (RN) derivative} $\phi$, i.e. \emph{the spectral density} of ${X}$, such that for any measurable subset $A$: $\mu(A)=\int_A\phi(\lambda)\,d\lambda$. 
If $H(X)=\overline{span}\; X\,\subset L^2(\Omega,\mathcal{F},\mathbf{P})$ the mean square estimation theory for stationary sequences is mainly based on \emph{Kolmogorov's isomorphism}:
\begin{equation}\label{iso}
I: L^2 (\mathbb{R},\mu)  \longrightarrow H({X})
\end{equation} 
given by the \emph {stochastic integral}:
$I(f)=\int\limits_{\mathbb{R}} f(\lambda)\,d\Phi(\lambda)\,,$ 
where $\Phi$ is the (orthogonal) random measure associated to ${X}$. 
In fact, if $A$ is a Borel subset then $\mu$ and $\Phi$ are related by the following formulas: 
$\mathbf{E}|I(\mathbf{1}_A)|^2 =\mathbf{E}|\Phi(A)|^2=\mu (A)$ and $\mathbf{E}\Big\vert{\int\limits_{\mathbb{R}}f\,d\Phi}\Big\vert^2=\int_{\mathbb{R}}|f|^2\,d\mu$. 
Moreover $X$ has the following representation:
\begin{equation}\label{repre} 
X(t)=I(e_t)=\int\limits_{\mathbb{R}} e^{it\lambda}\,d\Phi(\lambda)\,.
\end{equation}
with $e_t(\lambda)=e^{i t\lambda}$.
 
Analogous results hold for $X=(X(k))_{k\in\mathbb{K}}$ a w.s.s. random process with index time $k\in\mathbb{K}$ but with the integrals being taken over $\mathbb{T}:=\left[{0,\frac{2 \pi}{b}}\right)$, for example in this case $X(k)=\int\limits_{\mathbb{T}} e^{i k\lambda}\,d\Phi(\lambda)$. 
Indeed with appropriate restrictions the whole theory can be constructed for stationary processes indexed over more general Locally Compact Abelian (LCA) Groups. 
Following \cite{Roz}, linear time invariant filtering operations on ${X}$ are defined by:
\begin{equation}\label{filterfreq}
Y(t)=\int\limits_{\mathbb{R}} {f(\lambda)} e^{it\lambda}\,d\Phi(\lambda),\hspace{0.5cm}f\in L^2 (\mathbb{R}, d\mu)\,,
\end{equation}
so the resulting stationary process ${Y}=(Y(t))_{t\in\mathbb{R}}$ can be thought as the output of a linear system with a frequency response given by $f$ (i.e. \emph{filter}) and a random input $X$. 
In this case, the covariance of $Y$ is given by:
\begin{equation}\label{covfil}
R_Y(t-u)=\mathbf{E}(Y(t)\overline{Y(u)})=\int_{\mathbb{R}} |f(\lambda)|^2 e^{i\lambda(t-u)}\,d\mu (\lambda)\,.
\end{equation}
Finally, the spectral measure $\mu$ can be decomposed into a continuous and purely discrete part $\mu_c$ and $\mu_d$ and there exist measurable subsets $C,D$ such that $\mu_c(A)=\mu(A\cap C)$ and $\mu_d (A)=\mu(A\cap D)$. 
From this we can give an orthogonal (independent in the Gaussian case) decomposition of $X$, 
\begin{equation}\label{descDC}
X(t)=\int\limits_{C} e^{it\lambda} d\Phi(\lambda) + \int\limits_{D} e^{it\lambda}\,d\Phi(\lambda)=X_c (t) + X_d(t)\hspace{0.5cm}\textnormal{a.s.}\,. 
\nonumber
\end{equation}
This corresponds to the case when one replaces $f=\mathbf{1}_C$, or $f=\mathbf{1}_D$, in \eqref{filterfreq}. 
For short, $X_c$ and $X_d$ will be called the \emph{continuous} and \emph{discrete} parts of $X$ respectively. 
If $\mu_X$ is discrete then it is concentrated over a (countable) subset $D_X$ of $\mathbb{R}$ and moreover \eqref{repre} takes the form of a random series:
\begin{equation}\label{repredisc}
X(t)=\sum\limits_{\lambda\in D_X} e^{it\lambda}\Phi(\{\lambda\})\,.
\end{equation}
For short we sometimes will write $C(\lambda)=\Phi({\lambda})$ for these random coefficients. 
In this case, we shall say that the process $X$ has \emph{discrete spectrum}, and in contrast if $\mu_X= \mu_{X\,c}$ we will say that $X$ has \emph{continuous spectrum}. 
 
\subsubsection{The Ergodic Theorems}\label{ergosec} 

Natural estimators of the mean, variance and other statistics of a stationary process $X$ are appropriate time averages. 
We say that the strictly stationary process $X$ is \emph{metrically transitive} if the only measurable sets which are invariant under the \emph{shift} $X \longmapsto X^T=(X(t+T))_{t\in\mathbb{R}}$ have probability zero or one. 
If we consider the measurable space $(\mathbb{R}^{\mathbb{R}},\Sigma)$ where $\Sigma$ is the $\sigma$-algebra generated by the cylinder sets, then Birkhoff's ergodic theorem states (See e.g. \cite[p. 76]{Dym} or \cite[Theorem 5.1, p. 157]{Roz}):

\begin{thm}\label{ergo1} 
Let $X$ be a strictly stationary random process. 
If $f$ is a $(\mathbb{R}^{\mathbb{R}},\Sigma)$-measurable function and $\mathbf{E}|f(X)|<\infty$ then:
$$\lim\limits_{T\longrightarrow\infty}\frac{1}{T}\,\int_0 ^T f(X^t)\,dt=\mathbf{E}(f(X)|\mathcal{S})\hspace{0.5cm}\textnormal{a.s.}$$ 
and in the norm of $L^p(\Omega,\mathcal{F},\mathbf{P})$ for $1\leq p<\infty$, where $\mathcal{S}$ is the $\sigma$-algebra of shift-invariant sets.
\end{thm}
One can deduce that $X$ is metrically transitive if and only if in the above limit $\mathbf{E}(f(X)|\mathcal{S})=\mathbf{E}(f(X))$ for every such $f$. 
Metric transitivity, if $X$ is Gaussian, can be elegantly described in terms of its spectral measure $\mu$ by Maruyama's theorem (See e.g. \cite[p. 76]{Dym}  or \cite[p.163-166]{Roz}):

\begin{thm}\label{maru} 
A stationary Gaussian process $X$ is metrically transitive if and only if 
its spectral measure $\mu$ is continuous.
\end{thm}
If $X$ is only wide sense stationary, Von Neumann's mean ergodic theorem 
gives also an answer (See e.g. \cite[Theorem 6.2, p. 25]{Roz}):

\begin{thm}\label{VNeu} 
If $X$ is a w.s.s. random process then the equality
$$\lim\limits_{T\longrightarrow\infty}\frac{1}{T}\,\int_0 ^T X(t)\,dt=\mathbf{E}(X(0))$$
holds if and only if $\mu(\{0\})=0$. 
The limit is taken in the mean square sense.
\end{thm}
These results can be adapted to the two-sided averages since 
$$\int_{-T}^T f(X^t)\,dt=\int_0 ^T f(X^t)\,dt + \int_0 ^T f(X^{-t})\,dt\,.$$ 
We shall apply directly this device without any mention of it. 
In particular we will consider the averages given by the Besicovitch's norms $\lim\limits_{T\longrightarrow\infty}\frac{1}{2T}\int_{-T}^{T}|X(t)|^2\,dt=\Vert X\Vert^2_{B^2(\mathbb{R})}$ or $\Vert X\Vert^2_{B^2(\mathbb{K})}$ in the discrete case, since the same results hold for stationary random sequences replacing the integrals with sums. 
Indeed, as a consequence of this fact we have also the following useful result:

\begin{prop}\label{pro1}
Let $X$ be a zero mean Gaussian stationary random processes with continuous spectrum. 
Then for all $\tau\in \mathbb{R}$ the limit $$\rho_X (\tau)=\lim\limits_{T\longrightarrow
\infty} \frac{1}{2T}\int_{-T}^{T} X(t) \overline{X(t+\tau)}\,dt$$ exists \textnormal{a.s.} and equals $\mathbf{E}(X(0)\overline{X(\tau)})=R_X (\tau)$.
\end{prop}
For a proof and additional details, see a more general result for LCA groups in \cite{Med}.

\subsubsection{Stationary processes and the space $B^2(\mathbb{R})$}\label{Sec-AP-frame}

The space $B^2 (\mathbb{R})$ is related to stationary random processes. In this context, we can take an apparent advantage if $X$ is a w.s.s. mean square continuous random process since we can consider a measurable and equivalent process \cite{Gikh}. 
If additionally $X$ is stationary, note that by Theorem \ref{ergo1}, the value of $\Vert X\Vert^2_{B^2(\mathbb{R})}$ (or $\Vert X\Vert^2_{B^2(\mathbb{K})}$ in the case of a random sequence) always exists and is finite a.s.. 
Wide sense stationary random processes are always the mean square limit of random trigonometric polynomials.
However we note that $X$ cannot always be described as an element of $B^2(\mathbb{R})$ or equivalently it is not the limit in $B^2(\mathbb{R})$-norm of trigonometric polynomials. 
In fact one has the following claim. 
The proof of this can be found in \cite{Med} in the context of LCA groups.
   
\begin{lem}\label{struct} 
Let $X$ be a zero mean Gaussian stationary random process. 
Then:\hfill\\
(i) If $X$ has discrete spectrum then $X\in\,B^2 (\mathbb{R})$ \textnormal{a.s.}. 
Moreover, if \eqref{repredisc} is the spectral representation of $X$, then: \eqref{repredisc} converges \textnormal{a.s.} for all $t\in\,\mathbb{R}$ and $\left\|\,X\,\right\|_{B^2 (\mathbb{R})}=\sum\limits_{\lambda\in \Lambda} | C (\lambda)|^2 $ \textnormal{a.s.}. \hfill\\
(ii) 
If $X$ has continuous spectrum and $X\in\, B^2 (\mathbb{R})$ \textnormal{a.s.} then $X$ is the trivial null process, i.e. for every $t\in\,\mathbb{R}$, $X(t)=0$ \textnormal{a.s.}.
\end{lem}

\subsubsection{Stochastic Integration}

Finally, we state a result on the interchange of an stochastic integral with an ordinary integral with respect to the Lebesgue measure. 
This simple lemma is an adaptation of one presented in \cite{Gikh}, so the proof is left to the reader. 
The measures $\Phi$ and $\mu$ are as described before.
\begin{lem}\label{fubstoch} \cite[p. 237]{Gikh}. 
Let $g(\lambda,t)$ and $h(t)$ be two Borel measurable functions such that:
\begin{equation}\label{cond1}
\int_{\mathbb{R}} \int_{\mathbb{R}} |g(\lambda,t)|^2\,d\mu(\lambda)\,dt<\infty\hspace{0.5cm}\textnormal{and}\hspace{0.5cm}h\in\,L^2(\mathbb{R})\,,
\nonumber
\end{equation}
then:
$$\int_{\mathbb{R}} h(t)\Bigg(\,{\int\limits_{\mathbb{R}}g(\lambda,t)\,d\Phi(\lambda)}\Bigg)\,dt=
\int\limits_{\mathbb{R}} \left({\int_{\mathbb{R}} h(t)g(\lambda,t)\,dt}\right)\,d\Phi(\lambda)\hspace{0.5cm}\textnormal{a.s.}\,.$$
\end{lem}
The following examples of application of this result will be useful in the future. 
For a detailed justification see e.g. \cite{CenMed22}.
 
\subsubsection{Examples}\label{ex1} 
\begin{itemize}
\item Let $X$ be a w.s.s. process with a spectral representation given by equation (\ref{repre}). 
If $f\in L^1(\mathbb{R})$, it is possible to write an alternative expression for the inner product $\langle X,f\rangle$, which now defines a random variable. 
Indeed,
\begin{equation}\label{fourstoch}
\langle X,f\rangle=
\int\limits_{\mathbb{R}}\overline{\widehat{f}}(\lambda)\,d\Phi(\lambda)\hspace{0.5cm}\textnormal{a.s.}\,.
\end{equation}
\item From equation \eqref{fourstoch}, we can give a more intuitive interpretation of the filtering operation of equation \eqref{filterfreq} as a convolution. 
Suppose that $f=\widehat{\varphi}$, with $\varphi\in L^1 (\mathbb{R})$. 
Then given $t\in\mathbb{R}$:
$$Y(t)=(X\ast\varphi)(t)=\int_{\mathbb{R}}X(s)\varphi(t-s)\,ds=\int\limits_{\mathbb{R}}{\widehat{\varphi}(\lambda)}e^{i\lambda t}\,d\Phi(\lambda)\,.$$
\item If in addition $X$ is Gaussian with discrete spectrum, that is \eqref{repre} takes the form of the series $X(t)= \sum_{\lambda\in\Lambda} C(\lambda) e^{i\lambda t}$ with $C(\lambda)$ uncorrelated (thus independent) Gaussian random variables, and in this case:
\begin{equation}\label{ex2}
\langle X,f\rangle= \int\limits_{\mathbb{R}}\overline{\widehat{f}}(\lambda)\,d\Phi(\lambda)= 
\sum_{\lambda\in\Lambda}\overline{\widehat{f}}(\lambda)
C (\lambda)\,,\end{equation}
where the right hand series converges a.s..
\end{itemize}
 
\section{AP-frames and stationary processes}\label{AP-frames-processes}

Along this section $\mathcal{A}\subset L^1(\mathbb{R})\bigcap L^2(\mathbb{R})$ will be an affine system which verifies the additional conditions:
\begin{align*}
(C1)&\hspace{2cm}M_{\psi,j}:=\sup_{\lambda\in\mathbb{T}}\sum_{d\in\mathbb{D}}|\widehat{\psi}(a^j(\lambda +d))|^2\leq M_\psi <\infty\hspace{0.5cm}\forall\,j\in\mathbb{Z}, \\ 
(C2)&\hspace{2cm} \mathds{G}(\cdot)(0,q)\in C(\mathbb{R}\smallsetminus \{0,q\})\;\textnormal{ for each }q\in Q. 
\end{align*}
Under these assumptions, given $X$ a Gaussian stationary random process, recalling Examples \ref{ex1} (equation \eqref{fourstoch}) we can calculate its \emph{random} frame coefficients and get, for each $j\in\mathbb{Z}$, the random sequence $(\langle{X,\psi_{j,k}}\rangle)_{k\in\mathbb{K}}$ given by:
\begin{equation}\label{coefs}
\langle{X,\psi_{j,k}}\rangle= \int\limits_{\mathbb{R}}e^{i\lambda ka^j}\,\overline{\widehat{\psi}}(a^j\lambda)\,d\Phi (\lambda),
\end{equation}
It is immediate that these sequences belong to $H(X)$ and are also Gaussian. 
Moreover, for each $j$, $(\langle{X,\psi_{j,k}}\rangle)_{k\in\mathbb{K}}$ is stationary, since from \eqref{covfil} the covariance of this sequence takes the form:
\begin{equation}\label{covcoef}
\mathbf{E}(\langle{X,\psi_{j,k}}\rangle \overline{\langle{X,\psi_{j,k'}}\rangle})=\int_{\mathbb{R}}e^{i\lambda (k-k')a^j}|\widehat{\psi}(a^j\lambda)|^2\,d\mu(\lambda),
\end{equation}
which only depends on the difference $k-k'\in\mathbb{K}$. 
We shall characterize $\mathcal{A}$ with these coefficients. 
The apriori continuity conditions like (C2), as noted in \cite{Kim2009} for the case of $AP(\mathbb{R})$, are due to the fact that the $L^2(\mathbb{R})$ setup is associated with the Lebesgue measure in the frequency domain, while in the present setup the frequency domain is associated to stationary random processes with an arbitrary finite Borel spectral measure $\mu$. 
As a result, conditions like \eqref{RSpos} which are valid a.e. with respect to the Lebesgue measure will be valid \emph{for every} $\lambda$. 
We recall from \cite{CenMed22} this auxiliary lemma before proving the main results.

\begin{lem}\label{bobo} 
Let $f\in C(\mathbb{R}^N)$, $X:\Omega\longrightarrow\mathbb{R}^N$ be a random vector with a joint probability density $p_X$ such that $supp(p_X)=\mathbb{R}^N$ and $C\subseteq\mathbb{R}$ a closed subset such that $\mathbf{P}(f(X)\in C)=1$ then $f(x)\in C$ for all $x\in\mathbb{R}^N$. 
\end{lem}
Now, we can prove:

\begin{lem}\label{l1} 
Let $\mathcal{A}$ be an affine system. 
Then $\mathcal{A}$ is an $L^2(\mathbb{R})$-frame if and only if there exist constants $0<A\leq B$ such that:
\begin{equation}\label{e1}
A\Vert X\Vert^2_{B^2(\mathbb{R})}\leq
\sum\limits_{j\in\mathbb{Z}}\lim_{N\longrightarrow\infty}\frac{1}{(2N+1)}\sum\limits_{k\in\mathbb{K}(N)}|\langle{X, \psi_{j,k}}\rangle|^2 \leq 
B\Vert X\Vert^2_{B^2(\mathbb{R})}\hspace{0.5cm}\textnormal{a.s.} .
\end{equation}
for every Gaussian stationary random process $X$ with discrete spectral measure $\mu$ satisfying $\mu(\{0\})=0$.
\end{lem}
\begin{proof} (``Only if'' part) 
Recall from the beginning of this section that if $(X(t))_{t\in\mathbb{R}}$ is stationary and Gaussian, then for each $j\in\mathbb{Z}$, the sequence $(\langle{X,\psi_{j,k}}\rangle)_{k\in\mathbb{K}}$, given by $\langle{X,\psi_{j,k}}\rangle= \int\limits_{\mathbb{R}}e^{i\lambda ka^j}\,\overline{\widehat{\psi}}(a^j\lambda)\,d\Phi (\lambda)$ is also stationary (see equation \eqref{covcoef}) and Gaussian. 
If $X$ has discrete spectrum there exists $\Lambda\subset\mathbb{R}\smallsetminus\{0\}$ countable or finite, such that 
$X(t)=\sum_{\lambda\in\Lambda}C(\lambda)e^{i\lambda t}$ converges a.s. for all $t\in\mathbb{R}$, where the $C(\lambda)$'s are normal, zero mean and independent random variables such that $\sum_{\lambda\in\Lambda}\mathbf{E}|C(\lambda)|^2 =\mu(\mathbb{R}) <\infty$. 
Moreover, by Lemma \ref{struct}, $X\in B^2(\mathbb{R})$ a.s..   
In this case, the $C(\lambda)=0$ a.s. except for finite or countable $\lambda$'s, then the coefficients $\langle{X,\psi_{j,k}}\rangle$ of equation \eqref{coefs} take the form of an a.s. convergent series for each $j$ and $k$:
\begin{equation}\label{coef1a}
\langle{X,\psi_{j,k}}\rangle=\sum_{\lambda\in\mathbb{R}}e^{i\lambda ka^j}\,\overline{\widehat{\psi}}(a^j\lambda)C (\lambda)\,,
\nonumber 
\end{equation}
since this is a particular case of  \eqref{ex2},
and the same argument holds for any rearrangement of this sum. 

For each $j\in\mathbb{Z}$ let 
$\mathbb{D}_j:=a^{-j}\mathbb{D}$ and $\mathbb{T}_j:=a^{-j}\mathbb{T}$. 
Then, for each $j\in\mathbb{Z}$ and $k\in\mathbb{K}$, by periodization, we have the a.s. equality
\begin{equation}\label{coef1}
\langle{X,\psi_{j,k}}\rangle=
\sum_{\lambda\in\mathbb{T}_j}\Bigg(\sum_{d\in\mathbb{D}_j}\overline{\widehat{\psi}}(a^j(\lambda+d))C(\lambda +d)\Bigg)e^{i\lambda ka^j}=\sum_{\lambda\in\mathbb{T}_j} D(\lambda,j) e^{i\lambda ka^j},
\end{equation}
since $e^{idka^j}=1$ for any $d\in\mathbb{D}_j$ and $k\in\mathbb{K}$. 
Therefore there exists $\Omega^s_j\in\mathcal{F}$ such that $\mathbf{P}(\Omega^s_j)=1$ where, for every $k\in\mathbb{K}$, the series (\ref{coef1}) converges. 
On the other hand, due to Theorem \ref{ergo1} there exists $\Omega^e_j\in\mathcal{F}$ such that $\mathbf{P}(\Omega^e_j)=1$ and where $\left\|{(\langle{X,\psi_{j,k}}\rangle)_{k\in\mathbb{K}}}\right\|_{B^2(\mathbb{K})}$ exists and is finite. 
Define for each $j\in\mathbb{Z}$: $\Omega_j =\Omega^e _j \cap\Omega^s _j$. 
Now, 
for each $j\in\mathbb{Z}$ the system $\big\{(e^{i\lambda ka^j})_{k\in\mathbb{K}}:\lambda\in\mathbb{T}_j\big\} $ is an orthonormal basis for $B^2(\mathbb{K})$. 
Thus equation (\ref{coef1}) implies that $D(\lambda,j)$ is the Bohr Transform of the sequence $(\langle{X,\psi_{j,k}}\rangle)_{k\in\mathbb{K}}\in B^2(\mathbb{K})$ over $\Omega_j$ and then Wiener's formula (\ref{wien2}) asserts that:
\begin{equation}
\left\|{(\langle{X,\psi_{j,k}}\rangle)_{k\in\mathbb{K}}}\right\|^2 _{B^2(\mathbb{K})}=
\left\|D(\cdot,j)\right\|^2 _{L^2(\mathbb{T}_j, dc)}=
\sum\limits_{\lambda\in\mathbb{T}_j}\vert D(\lambda,j)\vert^2\,,
\nonumber
\end{equation}
and then over $\bigcap_j\Omega_j$:
\begin{equation}
\sum\limits_{j\in\mathbb{Z}}\Vert{(\langle{X,\psi_{j,k}}\rangle)_{k\in\mathbb{K}}}\Vert^2 _{B^2(\mathbb{K})}=
\sum\limits_{j\in\mathbb{Z}}\sum\limits_{\lambda\in\mathbb{T}_j}\vert D(\lambda,j)\vert^2\,.
\nonumber
\end{equation}
We want to change the order of summation in the last equality. 
For this we will change the indexes of summation. 
To begin with we have
\begin{align*}
\sum\limits_{\lambda\in\mathbb{T}_j}\vert D(\lambda,j)\vert^2 &=
\sum\limits_{\lambda\in\mathbb{T}_j}\Bigg\vert \sum_{d\in\mathbb{D}_j}\overline{\widehat{\psi}}(a^j(\lambda+d))C(\lambda +d)
\Bigg\vert^2 \\&=
\sum\limits_{\lambda\in\mathbb{T}_j}\sum_{d,d'\in\mathbb{D}_j}C(\lambda +d)\overline{C}(\lambda +d')\overline{\widehat{\psi}}(a^j(\lambda+d))\widehat{\psi}(a^j(\lambda+d')).
\end{align*}
The boundedness assumption (C1) yields the estimation
\begin{align*}
\sum\limits_{\lambda\in\mathbb{T}_j}\sum_{d,d'\in\mathbb{D}_j}\vert C(\lambda +d)\overline{C}(\lambda &+d')\overline{\widehat{\psi}}(a^j(\lambda+d))\widehat{\psi}(a^j(\lambda+d'))\vert = \\&=
\sum\limits_{\lambda\in\mathbb{T}_j}\Bigg(\sum_{d\in\mathbb{D}_j}\vert C(\lambda +d)\overline{\widehat{\psi}}(a^j(\lambda+d))\vert\Bigg)^2 \\&\leq 
M_\psi \sum\limits_{\lambda\in\mathbb{T}_j}\sum_{d\in\mathbb{D}_j}\vert C(\lambda +d)\vert^2 =M_\psi\Vert X\Vert^2_{B^2(\mathbb{R})}<\infty \,,
\end{align*}
thus the series is absolutely (and unconditionally) convergent and hence, the following change in the summation order is valid
\begin{align*}
\sum\limits_{\lambda\in\mathbb{T}_j}\sum_{d,d'\in\mathbb{D}_j}C(&\lambda +d)\overline{C}(\lambda +d')\overline{\widehat{\psi}}(a^j(\lambda+d))\widehat{\psi}(a^j(\lambda+d'))=\\&= 
\sum_{d,d'\in\mathbb{D}_j}\sum\limits_{\lambda\in\mathbb{T}_j}C(\lambda +d)\overline{C}(\lambda +d')\overline{\widehat{\psi}}(a^j(\lambda+d))\widehat{\psi}(a^j(\lambda+d'))\,.
\end{align*}
We define the following equivalence relation on $\mathbb{R}$: $\lambda\sim_Q\lambda'$ if and only if $\lambda -\lambda'\in Q:=\cup_j\mathbb{D}_j$. 
Let $\Delta$ be a system of representatives of $\sim_Q$. 
Thus, if $\lambda\in\mathbb{R}$ then there exist unique $\delta\in\Delta$ and $q\in Q$ such that $\lambda=\delta+q$. 
By straightforward calculations, for each $j\in\mathbb{Z}$ the following sets are identical. 
\begin{align*}
R_j&:=\{(x,y)\in\mathbb{R}^2:x=\lambda +d ,y=\lambda+d',\lambda\in\mathbb{T}_j,\,d,d'\in\mathbb{D}_j\},\\
S_j&:=\{(x,y)\in\mathbb{R}^2:x=\delta +q ,y=\delta +q',\delta\in\Delta,\,q-q'\in\mathbb{D}_j\}.
\end{align*}
These yield to
\begin{align*}
\sum\limits_{\lambda\in\mathbb{T}_j}\vert D&(\lambda,j)\vert^2 =\\&=
\sum\limits_{\lambda\in\mathbb{T}_j}\sum_{d,d'\in\mathbb{D}_j}C(\lambda +d)\overline{C}(\lambda +d')\overline{\widehat{\psi}}(a^j(\lambda+d))\widehat{\psi}(a^j(\lambda+d'))\\&=
\sum\limits_{\delta\in\Delta}\sum_{q,q'\in Q}C(\delta+q)\overline{C}(\delta+q')\overline{\widehat{\psi}}(a^j(\delta+q))\widehat{\psi}(a^j(\delta+q'))\mathbf{1}_{\mathbb{D}_j}(q-q')\\&=
\sum\limits_{\delta\in\Delta}\langle\mathds{G}_j(\delta)C_\delta, C_\delta\rangle_{\ell^2(Q)},
\end{align*}
where $C_\delta :=(C(\delta + q))_{q\in Q}\in\ell^2(Q)$. 

Now, by summing over $j\in\mathbb{Z}$, since $\langle\mathds{G}_j(\delta)C_\delta, C_\delta\rangle_{\ell^2(Q)}\geq 0$ for each $j\in\mathbb{Z}$, we can change the order of summation (by Fubini's Theorem) to obtain
\begin{align*}
\sum_{j\in\mathbb{Z}}\Vert(\langle X,\psi_{j,k}\rangle)_{k\in\mathbb{K}}\Vert^2_{AP(\mathbb{K})} &=
\sum\limits_{j\in\mathbb{Z}}\sum\limits_{\lambda\in\mathbb{T}_j}\vert D(\lambda,j)\vert^2 \\&=
\sum\limits_{j\in\mathbb{Z}}\sum\limits_{\delta\in\Delta}\langle\mathds{G}_j(\delta)C_\delta, C_\delta\rangle_{\ell^2(Q)}\\&=
\sum\limits_{\delta\in\Delta}\sum\limits_{j\in\mathbb{Z}}\langle\mathds{G}_j(\delta)C_\delta, C_\delta\rangle_{\ell^2(Q)}\\&=
\sum\limits_{\delta\in\Delta}\langle\mathds{G}(\delta)C_\delta, C_\delta\rangle_{\ell^2(Q)},
\end{align*}
where in the last equality is used \eqref{sumagrammianitos}.
Since we are assuming that $\mathcal{A}$ is an $L^2(\mathbb{R})$-frame, by applying Theorem \ref{RS} we have
\begin{equation}
A\Vert C_\delta\Vert^2_{\ell^2(Q)} \leq \langle\mathds{G}(\delta)C_\delta, C_\delta\rangle_{\ell^2(Q)}\leq B\Vert C_\delta\Vert^2_{\ell^2(Q)}.
\nonumber
\end{equation}
The result follows by noting that
\begin{equation}
\sum_{\delta\in\Delta}\Vert C_\delta\Vert^2_{\ell^2(Q)}=\sum_{\delta\in\Delta}\sum_{q\in Q}\vert C(\delta+q)\vert^2 =\sum_{\lambda}\vert C(\lambda)\vert^2 = \Vert X\Vert^2_{B^2(\mathbb{R})}.
\nonumber
\end{equation}
(``If'' part) Let $\tilde{Q}$ be a finite subset of $Q$ and define the random process $X(t)=\sum_{q\in\tilde{Q}}C(q)e^{i(q+\lambda)t}$, with $C(q)\sim\mathcal{N}(0,\sigma_q ^2)$ independent random variables, $\lambda\in\mathbb{R}$. 
By a direct calculation, from the definition of Fourier Transform, for all $j\in\mathbb{Z}$ and $k\in\mathbb{K}$, we have:
\begin{equation}
\langle{X,\psi_{j,k}}\rangle=
\sum\limits_{q\in\tilde{Q}} C(q)\int_{\mathbb{R}} e^{i(q+\lambda)t}\overline{\psi_{j,k}(t)}\,dt=
\sum\limits_{q\in\tilde{Q}}C(q)\overline{\widehat{\psi}}(a^j(q+\lambda))e^{i(q+\lambda)k a^j}\,,
\nonumber
\end{equation}
thus
\begin{equation}
|\langle{X,\psi_{j,k}}\rangle|^2=
\sum\limits_{q,q'\in\tilde{Q}}C(q)\overline{C}(q')\widehat{\psi}(a^j(q'+\lambda))\,\overline{\widehat{\psi}}(a^j(q+\lambda))e^{i(q-q')ka^j}\,,
\nonumber
\end{equation}
and summing over $k$, then for all $N\in\mathbb{N}$:
\begin{equation}
\frac{1}{(2N+1)}\sum\limits_{k\in\mathbb{K}(N)}\vert\langle{X,\psi_{j,k}}\rangle\vert^2=
\nonumber
\end{equation}
\begin{equation}\label{e2}
=\sum\limits_{q,q'\in\tilde{Q}}C(q)\overline{C}(q')\widehat{\psi}(a^j(q'+\lambda))\,\overline{\widehat{\psi}}(a^j(q+\lambda))\,
\frac{1}{(2N+1)}\sum\limits_{k\in\mathbb{K}(N)}e^{i(q-q')ka^j},
\end{equation}
Since for any $\lambda,\lambda'\in\mathbb{R}$ holds
\begin{equation}\label{e2-lim}
\lim_{N\longrightarrow\infty}\frac{1}{(2N+1)}\sum\limits_{k\in\mathbb{K}(N)}e^{i(\lambda-\lambda')ka^j}=\mathbf{1}_{\mathbb{D}_j}(\lambda-\lambda'),
\end{equation}
then
\begin{align*}
\lim_{N\longrightarrow\infty}\frac{1}{(2N+1)}&\sum\limits_{k\in\mathbb{K}(N)}\vert\langle{X,\psi_{j,k}}\rangle\vert^2 =\\&=
\sum\limits_{q,q'\in\tilde{Q}} C(q)\overline{C}(q')\widehat{\psi}(a^j(q'+\lambda))\,\overline{\widehat{\psi}}(a^j(q+\lambda))\mathbf{1}_{\mathbb{D}_j}(q-q')\,.
\end{align*}
On the other hand:
\begin{equation}\label{e3}
\lim_{T\longrightarrow\infty}\frac{1}{2T}\int_{-T}^T |X(t)|^2\,dt =\sum_{q\in\tilde{Q}} |C(q)|^2\,.
\end{equation}
Combining equations (\ref{e1}),(\ref{e2}),(\ref{e2-lim}),(\ref{e3}) and then summing over $j\in\mathbb{Z}$ it follows a.s. that 
\begin{equation}\label{e4}
A\sum_{q\in\tilde{Q}} |C(q)|^2\leq
\sum\limits_{q,q'\in\tilde{Q}}C(d)\overline{C}(q')\mathds{G}(\lambda)(q,q')
\leq B \sum_{q\in\tilde{Q}} |C(q)|^2\,,
\nonumber
\end{equation}
or equivalently, given $\lambda\in\mathbb{R}$ and $C(q)\sim\mathcal{N}(0,\sigma^2_q)$ independent random variables, $q\in\tilde{Q}$, if we define the finite random vector $Z(q)=C(q)\mathbf{1}_{\tilde{Q}} (q)$ it holds that
$$A\left\|Z\right\|^2_{\ell^2(\tilde{Q})}\leq\langle{\mathds{G}(\lambda) Z,Z}\rangle_{\ell^2(\tilde{Q})}\leq B\left\|Z\right\|^2_{\ell^2(\tilde{Q})}\,.$$
Recalling Lemma \ref{bobo} with $p_Z$ the Gaussian probability density of $Z$, one gets that for all $\lambda\in\mathbb{R}$ and $x\in\mathbb{R}^{\tilde{Q}}$: 
\begin{equation}\label{contM}
A\left\|x\right\|^2 _{\ell^2(Q)}\leq
\langle{ \mathds{G}(\lambda) x,x}\rangle_{\ell^2(Q)}\leq 
B \left\|x\right\|^2 _{\ell^2(Q)}\,.
\end{equation}
for all finite sequences $x$. 
Thus $\mathds{G}(\lambda)$ is bounded with continuous inverse, since the subset of finite sequences is dense in $\ell^2(Q)$. 
Then \eqref{contM} holds for arbitrary $x\in\ell^2 (Q)$. 
And the claim is proved recalling Theorem \ref{RS} (see \eqref{RSpos}).
\end{proof}
\emph{Remark.}
The condition $\mu(\{0\})=0$ was imposed since sums involving $\widehat{\psi}(a^j\lambda)$ may be undefined or can diverge for $\lambda=0$. 
Similar conditions are used in \cite{Kim2009}, in the context of almost periodic functions. 
Note that the general Definition \ref{AS} of an affine system does not require $\widehat{\psi}(0)=0$ as in the usual way it does Wavelet theory.

\begin{lem}\label{l2} 
Let $\mathcal{A}$ be an affine system. 
There exist constants $0<A\leq B$ such that $\mathcal{A}$ verifies equation (\ref{e1}) for every Gaussian stationary random process $X$ with continuous spectrum 
if and only if:
\begin{equation}\label{e5}
A \leq \sum\limits_{j\in\mathbb{Z}}|\widehat{\psi}(a^j\lambda)|^2\leq B
\end{equation}
for all $\lambda\neq 0$. \hfill\\
\textnormal{It is worth noting in this case that} $\mu(\{0\})=0$ \textnormal{and} $\Vert X\Vert^2_{B^2(\mathbb{R})}=\mu(\mathbb{R})$. 
\end{lem}
\begin{proof}
Let $X=(X(t))_{t\in\mathbb{R}}$ a Gaussian stationary random process with continuous spectrum. 
Then, recalling equation \eqref{covcoef}, for each $j\in\mathbb{Z}$ the sequence of coefficients $(\langle{X,\psi_{j,k}}\rangle)_{k\in\mathbb{K}}$ of equation \eqref{coefs}, is also stationary and Gaussian. 
Noting that this sequence has continuous spectrum (for a proof including this assertion see Lemma \ref{l3}), then by Theorem \ref{maru}, 
\begin{align}\label{e6}
\lim_{N\rightarrow\infty}\frac{1}{(2N+1)}\sum\limits_{k\in\mathbb{K}(N)}|\langle{X,\psi_{j,k}}\rangle|^2 &=
\mathbf{E}|\langle{X,\psi_{j,0}}\rangle|^2 \nonumber\\&=
\int_{\mathbb{R}}|\widehat{\psi}(a^j\lambda)|^2\,d\mu(\lambda)\hspace{0.5cm}\textnormal{a.s.}\,.
\end{align} 
A similar argument holds for the original process $X$:
\begin{equation}\label{e7}
\Vert X\Vert_{B^2(\mathbb{R})}^2=\lim_{T\longrightarrow\infty}\frac{1}{2T}\int_{-T}^T|X(t)|^2\,dt=\mathbf{E}|X(0)|^2=\mu(\mathbb{R})\hspace{0.5cm}\textnormal{a.s.}.
\end{equation} 
Now, set $\Omega_c$ as the subset where equation (\ref{e7}) holds, and for each $j$, $\Omega_j$ as the subset where equation (\ref{e6}) holds. 
Clearly, if $\mathop{\Omega}\limits^{\sim}:=\bigcap_{j\in\mathbb{Z}}(\Omega_c\cap \Omega_j)$ then $\mathbf{P}(\mathop{\Omega}\limits^{\sim})=1$ and consequently equation (\ref{e7}) is true over $\mathop{\Omega}\limits^{\sim}$ as well as 
\begin{align}\label{e8}
\sum\limits_{j\in \mathbb{Z}}\lim_{N\rightarrow\infty}\frac{1}{(2N+1)}\sum\limits_{k\in\mathbb{K}(N)}|\langle{X,\psi_{j,k}}\rangle|^2 &=
\sum\limits_{j\in\mathbb{Z}}\mathbf{E}|\langle{X,\psi_{j,0}}\rangle|^2 \nonumber\\ &=
\int_{\mathbb{R}}\sum\limits_{j\in \mathbb{Z}}|\widehat{\psi}(a^j\lambda)|^2\,d\mu(\lambda)\,.
\end{align}
Now, we can prove both implications. \hfill\\
If equation (\ref{e5}) is true, as equations (\ref{e7}) and (\ref{e8}) hold over $\mathop{\Omega}\limits^{\sim}$ is immediate that equation \eqref{e1} holds a.s..\hfill\\
Conversely, let $U \subset \mathbb{R}$ be any Borel subset of finite Lebesgue measure. 
Define the process $X(t):=\int\limits_U e^{i\lambda t}\,d\Phi(\lambda)$, where $\Phi$ is chosen as the Wiener random measure (see e.g. \cite[Example 2.3.1]{CenMed22} or p.13 of \cite{Roz}). 
Assuming that equation (\ref{e1}) holds, and noting that in this case $\mu$ is the restriction of the Lebesgue measure $m$ over $U$, equation (\ref{e7}) combined with (\ref{e8}) implies that:
$$A\,m(U) \leq  \int_U \sum\limits_{j\in\mathbb{Z}}|\widehat{\psi}(a^j\lambda)|^2\,d\lambda \leq B \, m(U)\,,$$
and then equation \eqref{e5} holds $m$-a.e., but the continuity condition (C2) proves the claim for every $\lambda$.  
\end{proof}

\emph{Remark.}
The careful reader can notice that one can state Lemma \ref{l2} for almost every $\lambda$ without condition (C2). 
However, this lemma will be used later combined with (C2) when dealing with more general random processes. \hfill\\ 
\hfill\\
We need the following lemmas from \cite{CenMed22} or \cite{Med} before we can prove one of the main results (Theorem \ref{t1}).

\begin{lem}\label{ortho}
i) Let $(Z(k))_{k\in \mathbb{K}}$ be a zero mean, Gaussian, stationary 
sequence. 
Then:
\begin{equation}\label{e9} 
\left\|Z\right\|^2_{B^2(\mathbb{K})}=\left\|Z_c\right\|^2_{B^2(\mathbb{K})} +\left\|Z_d\right\|^2 _{B^2(\mathbb{K})} \hspace{0.5cm}\textnormal{a.s.}\,.
\nonumber
\end{equation}
ii) Let $(Z(t))_{t\in \mathbb{R}}$ be a zero mean, Gaussian, stationary process. 
Then:
$$\left\|Z\right\|^2_{B^2(\mathbb{R})}= \left\|Z_c\right\|^2_{B^2(\mathbb{R})} + \left\|Z_d\right\|^2_{B^2(\mathbb{R})} \hspace{0.5cm}\textnormal{a.s.}.$$
\end{lem}

Recall that given $j\in\mathbb{Z}$ the sequence of coefficients $\langle X,\psi_{j,k}\rangle$ is also stationary. 
Then we have: 

\begin{lem}\label{l3} 
For each $j\in \mathbb{Z}$, if $k\in\mathbb{K}$ define $Z^j (k)=\langle X,\psi_{j,k}\rangle$. 
Then: 
$$Z^j _{c} (k)=\langle X_{c},\psi_{j,k}\rangle \hspace{0.5cm}\textnormal{and}\hspace{0.5cm}Z^j _{d}(k)=\langle X_d,\psi_{j,k}\rangle .$$
\end{lem}
\begin{proof}
First, note that for each $j$:
$$\langle X,\psi_{j,k}\rangle=\langle X_c,\psi_{j,k}\rangle + \langle X_d,\psi_{j,k}\rangle\,.$$
Let $\gamma$ and $\rho$ denote the spectral measures of the sequences $(\langle X_c,\psi_{j,k}\rangle)_{k\in\mathbb{K}}$ and $(\langle X_d,\psi_{j,k}\rangle)_{k\in\mathbb{K}}$ respectively. 
Recall that $\gamma$ is the unique measure such that 
$$\mathbf{E}(\langle{X_c,\psi_{j,k}}\rangle \overline{\langle{X_c,\psi_{j,0}}\rangle})=\int_{\mathbb{T}_j} e^{i\lambda k}\,d\gamma(\lambda)\,,$$
for all $k\in\mathbb{K}$. 
But if $\mu_c$ is the continuous part of $\mu$, the spectral measure of $X$, a direct calculation shows that,
\begin{align*}
\mathbf{E}(\langle{X_c,\psi_{j,k}}\rangle \overline{\langle{X_c,\psi_{j,0}}\rangle}) &=
\int_{\mathbb{R}} e^{i\lambda ka^j} |\widehat{\psi}(a^j\lambda)|^2\,d\mu_c(\lambda) \\&=
\sum\limits_{d\in\mathbb{D}_j}\int_{\mathbb{T}_j+d} e^{i\lambda ka^j} |\widehat{\psi}(a^j\lambda)|^2\,d\mu_c (\lambda)\,.
\end{align*}
Therefore, by uniqueness of the Fourier Transform, $\gamma$ is given by:
$$\gamma(A)=\sum\limits_{d\in\mathbb{D}_j}\int_{A+d} |\widehat{\psi}(a^j\lambda)|^2\,d\mu_c (\lambda)\,,$$
for all $A\in\mathcal{B}(\mathbb{T}_j)$. 
Then, in particular, for any $\lambda_0\in\mathbb{T}_j$ and taking $A=\{\lambda_0\}$, $\gamma(\{\lambda_0\})=0$. 
Hence $\gamma$ is continuous. 
On the other hand the discrete part of $\mu$ can be written as $\mu_d= \sum\limits_{\lambda\in\Lambda} c_{\lambda} \delta_\lambda$, with $\Lambda$ a countable subset of $\mathbb{R}$, $\delta_{\lambda}$ the unit mass measure concentrated on $\lambda$, and $c_{\lambda}>0$. 
Then
\begin{align*}
\mathbf{E}(\langle{X_d,\psi_{j,k}}\rangle\overline{\langle{X_d,\psi_{j,0}}\rangle}) &=
\int_{\mathbb{R}} e^{i\lambda ka^j} |\widehat{\psi}(a^j\lambda)|^2\,d\mu_d(\lambda)\\&=
\sum\limits_{\lambda\in\Lambda}c_\lambda\int_{\mathbb{R}}e^{i\lambda ka^j}|\widehat{\psi}(a^j\xi)|^2\,d\delta_\lambda(\xi)\,.
\end{align*}
This shows that if $A\in\mathcal{B}(\mathbb{T}_j)$, then 
$$\rho(A)=\sum\limits_{\lambda\in\Lambda}\sum\limits_{d\in\mathbb{D}_j} c_{\lambda+d}\int_{A}|\widehat{\psi}(a^j\xi)|^2\,d\delta_{\lambda+d}(\xi)\,,$$
which is clearly discrete.
\end{proof}

\begin{thm}\label{t1} 
Let $\mathcal{A}$ be an affine system. 
Then $\mathcal{A}$ is a $L^2(\mathbb{R})$-frame if and only if there exist constants $0<A\leq B$ such that equation \eqref{e1} holds
for every Gaussian stationary random process $X$.
\end{thm}
\begin{proof} 
First, we prove the "only if" part. 
Noting that $X(t)$ can be decomposed in two orthogonal (or independent indeed) processes $X_c(t)$ and $X_d(t)$ such that $X(t)=X_c (t) + X_d (t)$ a.s. with continuous and discrete spectral measure respectively, we get that:
$$\langle X,\psi_{j,k}\rangle=\langle X_c,\psi_{j,k}\rangle + \langle X_d,\psi_{j,k}\rangle\,.$$
By Lemma \ref{ortho} (i) and Lemma \ref{l3} then, with probability one, for each $j\in\mathbb{Z}$:
\begin{equation}\label{e12}
\left\|{(\langle X  ,\psi_{j,k}\rangle)_{k\in\mathbb{K}}}\right\|^2_{B^2(\mathbb{K})} =
\left\|{(\langle X_c,\psi_{j,k}\rangle)_{k\in\mathbb{K}}}\right\|^2_{B^2(\mathbb{K})} +
\left\|{(\langle X_d,\psi_{j,k}\rangle)_{k\in\mathbb{K}}}\right\|^2_{B^2(\mathbb{K})}\,,
\nonumber
\end{equation}
and therefore, since $\mathbb{Z}$ is countable, 
$$\sum\limits_{j\in\mathbb{Z}}\left\|({\langle X,\psi_{j,k}\rangle)_{k\in\mathbb{K}}}\right\|^2_{B^2(\mathbb{K})}=$$
\begin{equation}\label{e12b}
=\sum\limits_{j\in\mathbb{Z}}\left\|{(\langle X_c,\psi_{j,k}\rangle)_{k\in\mathbb{K}}}\right\|^2_{B^2(\mathbb{K})} + 
\sum\limits_{j\in\mathbb{Z}}\left\|{(\langle X_d,\psi_{j,k}\rangle)_{k\in\mathbb{K}}}\right\|^2_{B^2(\mathbb{K})}
\hspace{0.5cm}\textnormal{a.s.}.
\end{equation}
Similarly, by Lemma \ref{ortho} (ii), 
\begin{equation}\label{e13}
\left\|X\right\|^2_{B^2(\mathbb{R})}=
\left\|X_c\right\|^2_{B^2(\mathbb{R})} +\left\|X_d\right\|^2_{B^2(\mathbb{R})}\hspace{0.5cm}\textnormal{a.s.}.
\end{equation}
If $\mathcal{A}$ is a $L^2(\mathbb{R})$-frame then by Lemma \ref{l1} the result holds for $X_d$. 
Additionally, by Theorem \ref{RS}, taking the canonical vector $x=e_0$, $e_0 (q)=\mathbf{1}_{\{0\}}(q)$ in equation (\ref{RSpos}) we get $A\leq\langle{\mathds{G}(\lambda)e_0, e_0}\rangle_{\ell^2(Q)}=\sum\limits_{j\in\mathbb{Z}} |\widehat{\psi}(a^j\lambda)|^2 \leq B$, with $\left\|e_0\right\|_{\ell^2(Q)}=1$ and therefore equation (\ref{e5}) holds. 
Therefore, by Lemma \ref{l2}, equation (\ref{e1}) also holds for $X_c$. 
Finally, as the claim holds for $X_c$ and $X_d$ the desired conclusion follows by recalling equations (\ref{e12b}) and (\ref{e13}). \hfill\\
Conversely, the "if part" is a direct consequence of Lemma \ref{l1} since equation (\ref{e1}) in particular holds for every process with discrete spectral measure. 
\end{proof}

\section{Fractional processes associated to stationary Gaussian processes and smoothness analysis}\label{fractional}

It is relevant in several applications to study which statistical information can be extracted from the discrete and random coefficients $\langle X,\psi_{j\,k}\rangle$ \cite{Mas2}.   
In this section we will apply the results of Section \ref{AP-frames-processes} to characterize \emph{smoothness} of Gaussian stationary process $X$. 
We shall prove that the wavelet coefficients give a good description of the behavior of such processes. 
In fact, we will see that an appropriately weighted version of \eqref{e1} does this job in analogous way to some characterizations of classical Sobolev spaces, see e.g. Chapter 6 of \cite{Meyer}. 

\subsection{Some results on the regularity of w.s.s. random processes}\label{reg}

Smoothness can be described in several ways, which under certain conditions, may be equivalent.
Generally, this is achieved by studying the increments $X(t+h)-X(t)$, the decay of the spectral measure or the existence of the (norm) derivatives. 
Previous to introduction of the machinery of Section \ref{AP-frames-processes} we will make a brief but more precise description of what is meant here by \emph{smoothness}. 
We will see that this shares several features with the classical description given by Stein for the Lebesgue and Sobolev spaces in \cite{Stein}. 
The next result establishes the equivalence between the ``{smoothness}'' of the process $X$ and an integrability condition satisfied by the covariance function $R_X$, relating in this way the theory of \textit{singular integrals} \cite{Mall, samko, Stein} with fractional derivatives of a Gaussian stationary random process. 

\begin{thm}\label{fractional-Gaussian2}
Let $X$ be a w.s.s. random process with spectral measure $\mu$ and $0<\alpha< 1$. 
Then the following are equivalent:\hfill\\
i) $\displaystyle{\int_\mathbb{R}\vert\lambda\vert^{2\alpha}\,d\mu(\lambda)<\infty}$.\hfill\\
ii) $\displaystyle{\int_\mathbb{R}\vert R_X(0)-R_X(h)\vert\,\frac{dh}{\vert h\vert^{1+2\alpha}}<\infty}$. 
\end{thm}
\begin{proof}
First of all, from the definition of $R_X$ and since $X$ is stationary, note that 
\begin{align*}
\mathbf{E}\vert X(h)-X(0)\vert^2 &=
2\vert R_X(0)-R_X(h)\vert
\end{align*}
yields to
\begin{align*}
\int_\mathbb{R}\vert R_X(0)-R_X(h)\vert\,\frac{dh}{\vert h\vert^{1+2\alpha}}<\infty
\hspace{0.15cm}\Longleftrightarrow\hspace{0.1cm}
\int_\mathbb{R}\mathbf{E}\vert X(h)-X(0)\vert^2\,\frac{dh}{\vert h\vert^{1+2\alpha}}<\infty\,.
\end{align*}
By Kolmogorov's isomorphism \eqref{iso} we have
\begin{equation}\label{rep-for}
\int_\mathbb{R}\mathbf{E}\vert X(h)-X(0)\vert^2\,\frac{dh}{\vert h\vert^{1+2\alpha}}=\int_\mathbb{R}\Big(\int_{\mathbb{R}}\vert e^{ih\lambda}-1\vert^2\,d\mu(\lambda)\Big)\,\frac{dh}{\vert h\vert^{1+2\alpha}}.
\end{equation}
We define 
\begin{equation}\label{malliavin-fun}
F_\alpha(\lambda)=\int_\mathbb{R}\vert e^{ih\lambda}-1\vert^2\,\frac{dh}{\vert h\vert^{1+2\alpha}},\hspace{0.5cm}\lambda\in\mathbb{R}.
\nonumber
\end{equation}
By reasoning as in \cite[p. 144]{Mall} or \cite[Proposition 4, p. 140]{Stein}, there exists a constant $C_\alpha>0$ such that $F_\alpha(\lambda)=C_\alpha\vert\lambda\vert^{2\alpha}$ for all $\lambda\in\mathbb{R}$.\hfill\\
Since
\begin{align*}
\int_\mathbb{R}\mathbf{E}\vert X(h)-X(0)\vert^2\,\frac{dh}{\vert h\vert^{1+2\alpha}}=
\int_\mathbb{R}F_\alpha(\lambda)\,d\mu(\lambda) =C_\alpha\int_\mathbb{R}\vert\lambda\vert^{2\alpha}\,d\mu(\lambda),
\end{align*}
the equivalence follows by changing the order of integration in \eqref{rep-for}.
\end{proof}

The case $\alpha=1$ shows that the (usual) derivative process $\frac{dX}{dt}$ exists in the $L^2(\Omega,\mathcal{F},\mathbf{P})$-sense and justifies the formal manipulations:
\begin{equation}
\frac{dX}{dt}=\lim_{h\to 0}\frac{X(t+h)-X(t)}{h}=\frac{d}{dt}\Bigg(\int\limits_\mathbb{R}e^{it\lambda}\,d\Phi(\lambda)\Bigg)=\int\limits_\mathbb{R}\frac{\partial(e^{it\lambda})}{\partial t}\,d\Phi(\lambda).
\nonumber
\end{equation} 
This result which is well documented in the literature, e.g. \cite{Gikh,Roz}, and shares some basic features of the more general statement of the following Theorem \ref{fractional-Gaussian3}.
In fact, we finish this section with an analogue characterization of the smoothness of the process related to singular integrals and fractional derivatives, which enable us to work with a wider range of $\alpha\in (0,2)$. 
By the $\alpha$-th fractional derivative $D^\alpha X$ of a Gaussian stationary random process $X$ we mean, at least formally, the process given by
\begin{equation}\label{D-alfa}
(D^\alpha X)(t)=\int\limits_\mathbb{R}\vert\lambda\vert^\alpha e^{it\lambda}\,d\Phi(\lambda).
\end{equation}
By \eqref{filterfreq} when $\lambda\mapsto\vert\lambda\vert^\alpha$ belongs to $L^2(\mathbb{R},d\mu)$ then $D^\alpha X$ is a well-defined stationary random process. 
The next theorem resembles a classical result in the Theory of singular integrals and fractional Sobolev spaces (see for instance \cite[Chapter V]{Stein}).
In particular, this shows that \eqref{D-alfa} can be also interpreted, in a more classic way, as an appropriate limit of the increments of $X$.
  
\begin{thm}\label{fractional-Gaussian3}
Let $X$ be a w.s.s. random process with spectral measure $\mu$ and $0<\alpha < 2$. 
Then the following are equivalent:\hfill\\
i) $\displaystyle{\int_\mathbb{R}\vert\lambda\vert^{2\alpha}\,d\mu(\lambda)<\infty}$.\hfill\\
ii) $\displaystyle{\int_\mathbb{R}\mathbf{E}\vert X(t+h)+X(t-h)-2X(t)\vert^2\,\frac{dh}{\vert h\vert^{1+2\alpha}}<\infty}$. \hfill\\
iii) $\displaystyle{\lim_{\varepsilon\to 0}\int_{\vert h\vert\geq\varepsilon}(X(t+h)-X(t))\,\frac{dh}{\vert h\vert^{1+\alpha}}}$ exists in the $L^2(\Omega,\mathcal{F},\mathbf{P})$-sense for each $t\in\mathbb{R}$.
\end{thm}
\begin{proof}
By the same procedure as Theorem \ref{fractional-Gaussian2}, 
for any $t\in\mathbb{R}$ we have
\begin{equation}
\mathbf{E}\vert X(t+h)+X(t-h)-2X(t)\vert^2 = \int_\mathbb{R}\vert e^{ith}+e^{-ith}-2\vert^2\,d\mu(\lambda).
\nonumber
\end{equation}
$i)\hspace{0.15cm}\Leftrightarrow\hspace{0.15cm}ii)$ The proof is similar to that one of Theorem \ref{fractional-Gaussian2} and it is omitted.\hfill\\
$i)\hspace{0.15cm}\Rightarrow\hspace{0.15cm}iii)$ The hypothesis guarantees the well-definiteness of the process
\begin{equation}
Y_\alpha(t)=c(\alpha)\int\limits_\mathbb{R}\vert\lambda\vert^\alpha e^{i\lambda t}\,d\Phi(\lambda)=c(\alpha)(D^\alpha X)(t),
\nonumber
\end{equation}
where $c(\alpha)$ is a constant that will be specified further below. 
Also, $Y_\alpha$ is a mean square continuous \emph{w.s.s.} random process. 
For $0<\varepsilon<1$ consider the process 
\begin{equation}\label{regul-truc-hypersing-proc}
(D^\alpha_\varepsilon X)(t)=
\int_{\vert h\vert\geq\varepsilon}\frac{(X(t+h)-X(t))}{\vert h\vert^{1+\alpha}}\,dh.
\end{equation}
By doing the substitution $h=\varepsilon y$ we can rewrite \eqref{regul-truc-hypersing-proc} as
\begin{equation}\label{regul-truc-hypersing-proc2}
(D^\alpha_\varepsilon X)(t)=
\frac{1}{\varepsilon^\alpha}\int_{\vert y\vert\geq 1}\frac{(X(t+\varepsilon y)-X(t))}{\vert y\vert^{1+\alpha}}\,dy
\nonumber
\end{equation}
Now, by introducing the stochastic integral representation of $X$ followed by a change of order of integration (Lemma \ref{fubstoch}) and multiplying by $\frac{\vert\lambda\vert^\alpha}{\vert\lambda\vert^\alpha}\mathbf{1}_{\mathbb{R}\smallsetminus\{0\}}(\lambda)$, conveniently regrouped, yield to
\begin{equation}\label{regul-truc-hypersing-proc-SR}
(D^\alpha_\varepsilon X)(t)=
d(\alpha)\int\limits_\mathbb{R}\Big(\vert\lambda\vert^\alpha \widehat{K_\alpha}(\varepsilon\lambda)\Big)e^{it\lambda}\,d\Phi(\lambda),
\end{equation}
where 
\begin{equation}
\widehat{K_\alpha}(\lambda)=
\frac{1}{d(\alpha)\vert\lambda\vert^\alpha}\int_{\vert y\vert\geq 1}\frac{\big(e^{i\lambda y}-1\big)}{\vert y\vert^{1+\alpha}}\,dy ,\hspace{0.5cm}\lambda\in\mathbb{R}\smallsetminus\{0\},
\nonumber
\end{equation}
with $K_\alpha\in L^1(\mathbb{R},dt)$ and $d(\alpha)$ is a constant so that $\widehat{K_\alpha}(0)=\int_\mathbb{R}K_\alpha(t)\,dt=1$ \cite[Corollary to Lemma 3.16, Corollary to Theorem 3.19 \& Remark 3.20]{samko} (here is used $0<\alpha<2$). 
We declare $c(\alpha)=d(\alpha)$.  
Since $\widehat{K_\alpha}$ is a continuous function on $\mathbb{R}$ vanishing at infinity, it is bounded and satisfies $\lim\limits_{\varepsilon\to 0}\widehat{K_\alpha}(\varepsilon\lambda)=\widehat{K_\alpha}(0)=1$ for any $\lambda\in\mathbb{R}$. 
Thus, by the equality
\begin{align*}
\mathbf{E}\vert(D^\alpha_\varepsilon X)(t)-Y_\alpha(t)\vert^2 &=
c(\alpha)^2\,\mathbf{E}\Bigg\vert\int\limits_\mathbb{R}\Big(\vert\lambda\vert^\alpha\big[\widehat{K_\alpha}§(\varepsilon\lambda)-1\big]\Big)e^{it\lambda}\,d\Phi(\lambda)\Bigg\vert^2\\&=
c(\alpha)^2\int_\mathbb{R}\vert\lambda\vert^{2\alpha}\,\vert\widehat{K_\alpha}(\varepsilon\lambda)-1\vert^{2}\,d\mu(\lambda),
\end{align*}
we can apply Lebesgue's dominated convergence Theorem (recall that $\mu$ is a finite measure) to obtain the desired result. \hfill\\
$iii)\hspace{0.15cm}\Rightarrow\hspace{0.15cm}i)$ 
Suppose now that $\lim\limits_{\varepsilon\to 0}(D^\alpha_\varepsilon X)(t)$ exists in the $L^2(\Omega,\mathcal{F},\mathbf{P})$-sense for each $t\in\mathbb{R}$, where $(D^\alpha_\varepsilon X)(t)$ is as in \eqref{regul-truc-hypersing-proc-SR}. 
Thus, by Fatou's Lemma
\begin{align*}
\infty>\lim\limits_{\varepsilon\to 0}\mathbf{E}\vert(D^\alpha_\varepsilon X)(t)\vert^2 &\geq 
\int_\mathbb{R}\vert\lambda\vert^{2\alpha}\Big(\lim_{\varepsilon\to 0}\big\vert \widehat{K_\alpha}(\varepsilon\lambda)\big\vert^2\Big)\,d\mu(\lambda)\\&=
c(\alpha)^2\int_\mathbb{R}\vert\lambda\vert^{2\alpha}\,d\mu(\lambda)\,.\tag*{\qedhere}
\end{align*}
\end{proof}

\subsection{Regularity of Gaussian stationary random processes in terms of affine frames}\label{subsec-Regularity}

Given $\mathcal{A}:=\{\psi_{j,k}:j\in\mathbb{Z},\,k\in\mathbb{K}\}$ an affine system for $L^2(\mathbb{R})$ and $0<\alpha< 1$ we consider the system $I^\alpha\mathcal{A}:=\{I^\alpha\psi_{j,k}:j\in\mathbb{Z},\,k\in\mathbb{K}\}$, where $(I^\alpha \psi_{j,k})(t)$ denotes the \textit{Riesz potential} of $\psi_{j,k}$ of the order $\alpha$ on $t\in\mathbb{R}$ \cite[p. 37]{samko} given by the convolution operator
\begin{equation}
(I^\alpha\psi_{j,k})(t):=
\frac{1}{\gamma(\alpha)}\big(R_\alpha\ast\psi_{j,k}\big)(t),     
\nonumber
\end{equation}
with
\begin{equation}
\gamma(\alpha)=\pi^\frac{1}{2}2^\alpha\frac{\Gamma(\frac{\alpha}{2})}{\Gamma(\frac{1-\alpha}{2})}\hspace{0.5cm}\textnormal{and}\hspace{0.5cm}
R_\alpha (t)=\vert t\vert^{\alpha-1},\hspace{0.5cm}t\in\mathbb{R}\smallsetminus\{0\}.
\nonumber
\end{equation}
The constant $\gamma(\alpha)$ is chosen so that its Fourier transform is $\widehat{I^\alpha \psi_{j,k}}(\lambda)=\vert\lambda\vert^{-\alpha}\widehat{\psi_{j,k}}(\lambda)$ for $\lambda\in\mathbb{R}\smallsetminus\{0\}$ (see \cite[Lemma 1, p. 117]{Stein} or \cite[pp. 37-38]{samko}). 
  
Under the additional condition of being Gaussian, a equivalent characterization of the smoothness of $X$ (in any of the equivalent interpretations of the Subsection \ref{reg}) is given in the following two lemmas. 
This involves the decay of the frame coefficients $(\langle X,\psi_{j,k}\rangle)_{k\in\mathbb{K}}$. Recall that by the results of the previous section it is necessary and sufficient to study when $|\lambda|^{2\alpha}d\mu$ defines a finite measure.
 
\begin{lem}\label{fractional-discrete}
Let $X$ be a Gaussian stationary random process with discrete spectral measure $\mu$ satisfying $\mu(\{0\})=0$, suppose that $\mathcal{A}$ is an affine system such that $\mathcal{A}$ and $I^\alpha\mathcal{A}$ are $L^2(\mathbb{R})$-frames and $0<\alpha<1$. 
Then the following are equivalent:
\begin{enumerate}
\item[i)] $\displaystyle{\int_\mathbb{R}\vert \lambda\vert^{2\alpha}\,d\mu(\lambda)<\infty}$.
\item[ii)] $\displaystyle{\sum_{j\in\mathbb{Z}}a^{-2j\alpha}\Vert (\langle X,\psi_{j,k}\rangle)_{k\in\mathbb{K}}\Vert^2_{B^2(\mathbb{K})}<\infty \hspace{0.5cm}\textnormal{a.s.}}$.
\end{enumerate}
\end{lem}
\begin{proof}
We have seen in the proof of Lemma \ref{l1} that
\begin{align*}
\Vert(\langle X,\psi_{j,k}\rangle)_{k\in\mathbb{K}}\Vert^2_{B^2(\mathbb{K})} =
\sum\limits_{\delta\in\Delta}\langle\mathds{G}_j(\delta)C_\delta, C_\delta\rangle_{\ell^2(Q)},
\end{align*}
for each $j\in\mathbb{Z}$, thus
\begin{equation}\label{aux1}
\sum_{j\in\mathbb{Z}}a^{-2j\alpha}\Vert(\langle X,\psi_{j,k}\rangle)_{k\in\mathbb{K}}\Vert^2_{B^2(\mathbb{K})} =
\sum\limits_{\delta\in\Delta}\sum_{j\in\mathbb{Z}}a^{-2j\alpha}\langle\mathds{G}_j(\delta)C_\delta, C_\delta\rangle_{\ell^2(Q)}
\hspace{0.5cm}\textnormal{a.s.}.
\end{equation}
Now, 
\begin{align*}
&a^{-2j\alpha}\langle\mathds{G}_j(\delta)C_\delta, C_\delta \rangle_{\ell^2(Q)}=\\
&\hspace{0.5cm}=a^{-2j\alpha}\sum_{q,q'\in Q}C(\delta+q)\overline{C}(\delta+q')\overline{\widehat{\psi}}(a^j(\delta+q))\widehat{\psi}(a^j(\delta+q'))\mathbf{1}_{\mathbb{D}_j}(q-q')
\end{align*}
and
\begin{align*}
a^{-2j\alpha}& C(\delta+q)\overline{C}(\delta+q')\overline{\widehat{\psi}}(a^j(\delta+q))\widehat{\psi}(a^j(\delta+q'))=\\
&= C(\delta+q)\vert\delta+q\vert^\alpha\overline{\widehat{I^\alpha\psi}}(a^j(\delta+q))\overline{C}(\delta+q')\vert \delta+q'\vert^\alpha\widehat{I^\alpha\psi}(a^j(\delta+q'))\,.
\end{align*}
Define in first place:
\begin{align*}
D_\delta^\alpha(q):= D^\alpha(\delta+q)= C(\delta+q)\vert\delta+q\vert^\alpha,\hspace{0.5cm}\delta\in\Delta,q\in Q,
\end{align*}
note that these are precisely the coefficients of the process $D^\alpha X$ if this exists. Also write
\begin{align*}
\mathds{G}_j^\alpha(\delta)(q,q'):= \overline{\widehat{I^\alpha\psi}}(a^j(\delta+q))\widehat{I^\alpha\psi}(a^j(\delta+q'))\mathbf{1}_{\mathbb{D}_j}(q-q'),\hspace{0.3cm}\delta\in\Delta,q,q'\in Q,
\end{align*}
therefore we can rewrite equation \eqref{aux1} as

\begin{align*}
\sum_{j\in\mathbb{Z}}a^{-2j\alpha}\Vert(\langle X,\psi_{j,k}\rangle)_{k\in\mathbb{K}}\Vert^2_{B^2(\mathbb{K})}&=
\sum\limits_{\delta\in\Delta}\sum_{j\in\mathbb{Z}}\langle\mathds{G}_j^\alpha(\delta)D_\delta^\alpha, D_\delta^\alpha\rangle_{\ell^2(Q)}\\&=
\sum\limits_{\delta\in\Delta}\langle\mathds{G}^\alpha(\delta)D_\delta^\alpha,D_\delta^\alpha\rangle_{\ell^2(Q)}
\end{align*}
where the last equality follows by reasoning as in Lemma \ref{l1}. 
Since $I^\alpha\mathcal{A}$ is an $L^2(\mathbb{R})$-frame, there exist constants $0<A_\alpha\leq B_\alpha$ such that
\begin{equation}
A_\alpha \Vert D_\delta^\alpha\Vert^2_{\ell_2(Q)}\leq 
\langle\mathds{G}^\alpha(\delta)D_\delta^\alpha,D_\delta^\alpha\rangle_{\ell^2(Q)}
\leq B_\alpha \Vert D_\delta^\alpha\Vert^2_{\ell_2(Q)}.
\nonumber
\end{equation}
Now, the equality
\begin{equation}
\sum_{\delta\in\Delta}\Vert D_\delta^\alpha\Vert^2_{\ell_2(Q)}=\sum_{\delta\in\Delta}\sum_{q\in Q}\vert C(\delta+q)\vert^2\vert\delta+q\vert^{2\alpha}=\Vert D^\alpha X\Vert^2_{B^2(\mathbb{R})},
\nonumber
\end{equation}
yields to
\begin{align*}
\sum_{j\in\mathbb{Z}}a^{-2j\alpha}\Vert (\langle X,\psi_{j,k}\rangle)_{k\in\mathbb{K}}\Vert^2_{B^2(\mathbb{K})}<\infty \hspace{0.2cm}\textnormal{a.s.}
\hspace{0.4cm}\Longleftrightarrow\hspace{0.4cm}\Vert D^\alpha X\Vert^2_{B^2(\mathbb{R})}<\infty\hspace{0.2cm}\textnormal{a.s.}\,.
\end{align*}
Finally, we can prove the claim observing that if $X$ has discrete spectrum, then 
$(D^\alpha X)(t)$ takes the form of a series of Gaussian and independent random elements in $\mathcal{H}:=\overline{\textnormal{span}}\{e^{i\lambda t}:\lambda\in\Lambda\}\subset B^2(\mathbb{R})$ which is a separable Hilbert space. 
By noting that 
\begin{equation}
\mathbf{E}\Vert D^\alpha X\Vert^2_{B^2(\mathbb{R})}=\int_\mathbb{R}\vert \lambda\vert^{2\alpha}\,d\mu(\lambda),
\nonumber
\end{equation}
the implication $i)\hspace{0.15cm}\Leftrightarrow\hspace{0.15cm}ii)$ follows at once since: $D^\alpha X\in B^2(\mathbb{R})$ a.s. if and only if $\mathbf{E}\Vert D^\alpha X\Vert^2_{B^2(\mathbb{R})}<\infty$, since for Gaussian series of independent random elements in a separable Hilbert space, a.s. convergence and convergence in the $p$-mean ($p=2$ in this case) are equivalent (see e.g. \cite[p.56]{Kwap}). 
\end{proof}

If $X$ is a Gaussian stationary random process with continuous spectral measure, by the Ergodic Theorem we have
\begin{equation}\label{aux2}
\sum_{j\in\mathbb{Z}}a^{-2j\alpha}\Vert (\langle X,\psi_{j,k}\rangle)_{k\in\mathbb{K}}\Vert^2_{B^2(\mathbb{K})}= 
\sum_{j\in\mathbb{Z}}a^{-2j\alpha}\mathbf{E}\vert\langle X,\psi_{j,0}\rangle\vert^2\hspace{0.5cm}\textnormal{a.s.}
\end{equation}
The next result shows that Lemma \ref{fractional-discrete} holds true for Gaussian stationary random processes with continuous spectral measure.

\begin{lem}\label{fractional-continuous}
Let $X$ be a Gaussian stationary random process with continuous spectral measure $\mu$, suppose that $\mathcal{A}$ is an affine system such that $\mathcal{A}$ and $I^\alpha\mathcal{A}$ are $L^2(\mathbb{R})$-frames and $0<\alpha<1$. 
Then the following are equivalent:
\begin{enumerate}
\item[i)] $\displaystyle{\int_\mathbb{R}\vert \lambda\vert^{2\alpha}\,d\mu(\lambda)<\infty}$.
\item[ii)] $\displaystyle{\sum_{j\in\mathbb{Z}}a^{-2j\alpha}\Vert (\langle X,\psi_{j,k}\rangle)_{k\in\mathbb{K}}\Vert^2_{B^2(\mathbb{K})}<\infty \hspace{0.5cm}\textnormal{a.s.}}$.
\end{enumerate}
\end{lem}
\begin{proof}
To begin with, by \eqref{aux2} we have
\begin{equation}
\sum_{j\in\mathbb{Z}}a^{-2j\alpha}\Vert (\langle X,\psi_{j,k}\rangle)_{k\in\mathbb{K}}\Vert^2_{B^2(\mathbb{K})}=
\sum_{j\in\mathbb{Z}}a^{-2j\alpha}\mathbf{E}\vert\langle X,\psi_{j,0}\rangle\vert^2\hspace{0.5cm}\textnormal{a.s.}.
\nonumber
\end{equation}
Since
\begin{equation}
\mathbf{E}\vert\langle X,\psi_{j,0}\rangle\vert^2=\int_\mathbb{R}\vert\widehat{\psi}(a^j\lambda)\vert^2\,d\mu(\lambda),
\nonumber
\end{equation}
then 
\begin{align*}
\sum_{j\in\mathbb{Z}}a^{-2j\alpha}\Vert (\langle X,\psi_{j,k}\rangle)_{k\in\mathbb{K}}\Vert^2_{B^2(\mathbb{K})}&=
\sum_{j\in\mathbb{Z}}a^{-2j\alpha}\int_\mathbb{R}\vert\widehat{\psi}(a^j\lambda)\vert^2\,d\mu(\lambda)\\&=
\int_\mathbb{R}\vert\lambda\vert^{2\alpha}\sum_{j\in\mathbb{Z}}\dfrac{\vert\widehat{\psi}(a^j\lambda)\vert^2}{\vert a^j\lambda\vert^{2\alpha}}\,d\mu(\lambda)\hspace{0.5cm}\textnormal{a.s.}.
\end{align*}
Now, since $\widehat{I^\alpha\mathcal{A}}$ is an $L^2(\mathbb{R})$-frame, there exist constants $0<A_\alpha\leq B_\alpha$ such that
\begin{equation}
A_\alpha \leq \sum_{j\in\mathbb{Z}}\dfrac{\vert\widehat{\psi}(a^j\lambda)\vert^2}{\vert a^j\lambda\vert^{2\alpha}}\leq B_\alpha\hspace{0.5cm}\textnormal{a.e. }\lambda,
\nonumber
\end{equation}
thus \textnormal{a.s.} we have
\begin{displaymath}
A_\alpha\int_\mathbb{R}\vert\lambda\vert^{2\alpha}\,d\mu(\lambda)\leq
\sum_{j\in\mathbb{Z}}a^{-2j\alpha}\Vert (\langle X,\psi_{j,k}\rangle)_{k\in\mathbb{K}}\Vert^2_{B^2(\mathbb{K})}\leq
B_\alpha\int_\mathbb{R}\vert\lambda\vert^{2\alpha}\,d\mu(\lambda)\,,
\end{displaymath}
which proves the desired equivalence.
\end{proof}

Noting that $\mu$ is finite and combining the results of the previous section, we can give the following characterization in terms of the weight $(|\lambda|^2 +1)^{\alpha}$ in a closer form to the usual description of the Sobolev spaces in terms of the Fourier transform.

\begin{thm}\label{fractional-Gaussian}
Let $X$ be a Gaussian stationary random process with spectral measure $\mu$ satisfying $\mu(\{0\})=0$, suppose that $\mathcal{A}$ is an affine system such that $\mathcal{A}$ and $I^\alpha\mathcal{A}$ are $L^2(\mathbb{R})$-frames and $0<\alpha<1$. 
Then the following are equivalent:
\begin{enumerate}
\item[i)] $\displaystyle{\int_\mathbb{R}(\vert \lambda\vert^2+1)^{\alpha}\,d\mu(\lambda)<\infty}$.
\item[ii)] $\displaystyle{\sum_{j\in\mathbb{Z}}(a^{-2j}+1)^\alpha\Vert (\langle X,\psi_{j,k}\rangle)_{k\in\mathbb{K}}\Vert^2_{B^2(\mathbb{K})}<\infty \hspace{0.5cm}\textnormal{a.s.}}$.
\end{enumerate}
\end{thm}
\begin{proof}
Since $X$ is a Gaussian stationary random process and $\mathcal{A}$ is an $L^2(\mathbb{R})$-frame, Theorem \ref{t1} yields to
\begin{equation}\label{eq-L2frame}
\sum_{j\in\mathbb{Z}}\Vert (\langle X,\psi_{j,k}\rangle)_{k\in\mathbb{K}}\Vert^2_{B^2(\mathbb{K})}\leq B\Vert X\Vert^2_{B^2(\mathbb{R})}{}<\infty \hspace{0.5cm}\textnormal{a.s.}.
\end{equation}
$i)\hspace{0.15cm}\Rightarrow\hspace{0.15cm}ii)$ 
Since $\vert\lambda\vert^{2\alpha}\leq (\vert \lambda\vert^2+1)^{\alpha}$ for all $\lambda\in\mathbb{R}$, by hypothesis and monotony we obtain
\begin{equation}
\int_\mathbb{R}\vert\lambda\vert^{2\alpha}\,d\mu(\lambda)<\infty .
\nonumber
\end{equation}
Since $X$ can be decomposed as $X=X_d +X_c$, with $X_d$ and $X_c$ Gaussian stationary random process with discrete and continuous spectral measures $\mu_{X,d}$ and $\mu_{X,c}$ respectively, satisfying $\langle X,\psi_{j,k}\rangle=\langle X_d,\psi_{j,k}\rangle +\langle X_c,\psi_{j,k}\rangle$ for any $j\in\mathbb{Z}$, $k\in\mathbb{K}$ and also $\Vert(\langle X,\psi_{j,k}\rangle)_{k\in\mathbb{K}}\Vert^2_{B^2(\mathbb{K})}=\Vert(\langle X_d,\psi_{j,k}\rangle)_{k\in\mathbb{K}}\Vert^2_{B^2(\mathbb{K})} +\Vert(\langle X_c,\psi_{j,k}\rangle)_{k\in\mathbb{K}}\Vert^2_{B^2(\mathbb{K})}$ a.s. for any $j\in\mathbb{Z}$, then
\begin{equation}
\int_\mathbb{R}\vert\lambda\vert^{2\alpha}\,d\mu_{X.d}(\lambda)+\int_\mathbb{R}\vert\lambda\vert^{2\alpha}\,d\mu_{X,c}(\lambda)=
\int_\mathbb{R}\vert\lambda\vert^{2\alpha}\,d\mu(\lambda)<\infty .
\nonumber
\end{equation}
The result follows by applying Lemma \ref{fractional-discrete} and Lemma \ref{fractional-continuous} to $X_d$ and $X_c$ respectively together with \eqref{eq-L2frame}. 
Indeed, since $(a^{-2j}+1)^\alpha\leq 2^\alpha (a^{-2j\alpha}+1)$ for any $j\in\mathbb{Z}$, by monotony we obtain
\begin{align*}
\sum_{j\in\mathbb{Z}}(a^{-2j}+1)^\alpha &\Vert (\langle X,\psi_{j,k}\rangle)_{k\in\mathbb{K}}\Vert^2_{B^2(\mathbb{K})} \leq \\&\leq
2^\alpha\sum_{j\in\mathbb{Z}}a^{-2j\alpha}\Vert(\langle X,\psi_{j,k}\rangle)_{k\in\mathbb{K}}\Vert^2_{B^2(\mathbb{K})} \\&=
2^\alpha\sum_{j\in\mathbb{Z}}a^{-2j\alpha}\Vert(\langle X_d,\psi_{j,k}\rangle)_{k\in\mathbb{K}}\Vert^2_{B^2(\mathbb{K})} +\\&\hspace{0.5cm}+
2^\alpha\sum_{j\in\mathbb{Z}}a^{-2j\alpha}\Vert(\langle X_c,\psi_{j,k}\rangle)_{k\in\mathbb{K}}\Vert^2_{B^2(\mathbb{K})} <\infty
\hspace{0.5cm}\textnormal{a.s.}\,.
\end{align*}
$ii)\hspace{0.15cm}\Rightarrow\hspace{0.15cm}i)$ 
The reasoning is similar to the previous implication but reversing the order. 
First we use the inequality $a^{-2j\alpha}\leq (a^{-2j}+1)^\alpha$ for any $j\in\mathbb{Z}$ and then the inequality $(\vert\lambda\vert^2+1)^\alpha\leq 2^\alpha(\vert\lambda\vert^{2\alpha}+1)$ for all $\lambda\in\mathbb{R}$.
\end{proof}

Note that the utility of the last results depends on verifying previously that $I^\alpha\mathcal{A}$ is also a frame. 
We close this section with the following easy example which shows that this is not a great obstacle in principle: 

\begin{ex}
\textnormal{
Bandpass wavelets.
Recall that if ${\psi}$ is bandlimitted then the characterizations of Subsection \ref{subsec-Regularity} are reduced to the following: 
A necessary and sufficient condition for $\mathcal{A}$ to be an $L^2 (\mathbb{R})$-frame is that the following holds:
\begin{equation}\label{BLframe}
A \leq \sum_{j\in\mathbb{Z}}{\vert\widehat{\psi}(a^j\lambda)\vert^2} \leq B\,\;\;\;a.e.\,,
\end{equation}
where $A,B$ are the frame bounds of $\mathcal{A}$.
Let us assume additionally that $\psi$ is bandpass, i.e. there exists $\lambda_1 >\lambda_0 >0$ such that $\widehat{\psi}$ vanishes outside $D=\{ \lambda_0 <|\lambda|< \lambda_1\}$. 
This is the case of, for example, the Meyer wavelet, see e.g. Chapter 2 of \cite{Meyer}. 
Now we must verify \eqref{BLframe} for $\widehat{I^\alpha \psi_{j,k}}(\lambda)=\vert\lambda\vert^{-\alpha}\widehat{\psi_{j,k}}(\lambda)$, $\lambda\in\mathbb{R}\smallsetminus\{0\}$. 
In fact, we can estimate the lower frame bound as
\begin{equation}
A_\alpha=\frac{A}{\lambda_1 ^{2\alpha}} \leq  \sum_{j: \lambda_0 \leq |a^j\lambda|\leq\lambda_1 } \dfrac{\vert\widehat{\psi}(a^j\lambda)\vert^2}{\vert a^j\lambda\vert^{2\alpha}}=\sum_{j\in\mathbb{Z}}\dfrac{\vert\widehat{\psi}(a^j\lambda)\vert^2}{\vert a^j\lambda\vert^{2\alpha}} \hspace{0.5cm}\textnormal{a.e. }\lambda.
\nonumber
\end{equation}
The upper bound $B_{\alpha}=\frac{B}{\lambda_0 ^{2\alpha}}$ can be obtained similarly.}
\end{ex}

Finally, observe that the condition on the spectral measure $\mu (\{0\})=0$ was imposed in order to find conditions for a general affine system as introduced in Definition \ref{AS}, following the ideas of \cite{Kim2009}. 
However, in wavelet theory, is usual to assume that $\widehat{\psi}(0)=0$. 
This happens, for example, when the (discrete) affine system is obtained by sampling the corresponding integral version of the wavelet transform. 
In this case the condition of $\mu$ having no atom at $\lambda=0$ may be dropped in all the statements where this was used. 

\subsection*{Aknowledgements}
This work was funded by the  Universidad de Buenos Aires, Grant. No. 20020220400 \hfill
210BA and CONICET Grant PIP 2023 No. 11220220100077CO, Buenos Aires, Argentina. 

H.D. Centeno:\hfill\\
(https://orcid.org/0000-0003-2435-7170)\hfill\\
Universidad de Buenos Aires, Facultad de Ciencias Exactas y Naturales, Departamento de Matemática.\\
Universidad de Buenos Aires, Facultad de Ingeniería, Departamento de Matemática.\\
Av. Paseo Col\'on 850, 1er. piso, (C1063ACV), Buenos Aires CF, Argentina\\
Correspondence: hcenteno@fi.uba.ar

\hfill\\
J.M. Medina:\hfill\\
(https://orcid.org/0000-0003-0370-6837)\hfill\\
Universidad de Buenos Aires, Facultad de Ingeniería, Departamento de Matemática.\\ 
Pontificia Universidad Católica Argentina, Facultad de Ingeniería.\\
Instituto Argentino de Matem\'atica ``A. P. Calder\'on'' - CONICET.\\
Saavedra 15, 3er piso (1083),
Buenos Aires CF, Argentina\\
Correspondence: jmedina@conicet.gov.ar


\begin{thebibliography}{27}



\bibitem{CenMed22} Centeno H. D., Medina J. M.,  
{\textit{AP-frames and Stationary Random Processes}}, 
Appl. Comput. Harmon. Anal., 61:1-24, 2022.

\bibitem{Chris} Christensen O.,  
\emph{An Introduction to Frames and Riesz Bases}, 
2nd. Edition, Birkh\"{a}user, 2016.


\bibitem{Cor} Corduneanu C.,  
\emph{Almost Periodic Oscilations and Waves}, 
Springer, 2009.

\bibitem{Dym} Dym H., McKean H. P.,  
\emph{Gaussian Processes, Function Theory, and the Inverse Spectral Problem}, 
Dover, 2008.




\bibitem{galindo} Galindo F.,  
\emph{Some Remarks on: ``{On the windowed Fourier transform and wavelet transform of almost periodic functions}''}, 
[Appl. Comput. Harmon. Anal. 10 (2001), no. 1, 45–60; mr1808199] by J. R. Partington and B. \"{U}nalmı\c{s}. 
Appl. Comput. Harmon. Anal., 16(3):174–181, 2004.



\bibitem{Gikh} Gikhman I. I., Skorokhod A. V.,  
\emph{The Theory of Stochastic Processes}, 
Vol. I., Springer, Berlin, 2004.







\bibitem{Kantor} Kantorovich L. V. K., Akilov G. P.,  
\emph{Functional Analysis}, 
2nd. Edition, Pergamon Press, 1982.

\bibitem{Katz} Katznelson Y.,  
\emph{An Introduction to Harmonic Analysis}, 
Dover, 1976.

\bibitem{Kim2009} Kim Y. H., Ron A., 
\emph{Time Frequency Representations of Almost Periodic Functions}, 
Constr. Approx. 29, pp. 303-323, 2009.

\bibitem{Kim2013} Kim Y. H.,  
\emph{Representations of Almost-Periodic Functions Using Generalized Shift-Invariant Systems in $\mathbb{R}^d$}, 
J. Fourier Anal Appl. 19, pp. 857-876, 2013.

 
\bibitem{Kwap} Kwapie\'{n} S. , Woyczy\'{n}ski W.A.,
\emph{Random Series and Stochastic Integrals: Single and Multiple}, 
Birkh\"{a}user, 1992.


\bibitem{Lenz2020} {D. Lenz, T. Spindeler, N. Strungaru}, 
\emph{Pure Point Diffraction and Mean, Besicovitch and Weyl Almost Periodicity}, 
arXiv:2006.10821v1, https://arxiv.org/abs/2006.10821, 2020.



\bibitem{Mall} Malliavin P., 
\emph{Integration and Probability}, 
Springer New York, 1995.

\bibitem{Mas2} Masry E., Cambanis S.,  
\emph{The Representation of Stochastic Processes Without Loss of Information}, 
SIAM J. Appl. Math., Vol. 25, No. 4, pp. 628-633, 
1973.





\bibitem{Med} Medina J.M., 
\emph{Laws of Large Numbers, Spectral Translates and Sampling over LCA Groups}. J. Fourier Anal. Appl. 30 (67), pp.1-32. 2024.

\bibitem{Meyer} Meyer Y. 
\emph{Wavelets and Operators}. 
Cambridge Studies in Advanced Mathematics, 37, 1992.

\bibitem{Partin} Partington J. R., \"{U}nalmi\c{s} B.,  
\emph{On the Windowed Fourier Transform and Wavelet Transform of Almost Periodic Functions}, 
Appl. Comput. Anal. Appl. 10, pp. 45-60, 2001. 

\bibitem{Pinsky} Pinsky M., 
\emph{Introduction to Fourier Analysis and Wavelets (Graduate Studies in Mathematics)}, 
American Mathematical Society, 2009.

\bibitem{Pour} Pourahmadi M.,  
\emph{Foundations of Time Series Analysis and Prediction Theory}, 
Wiley series in Statistics, 2001.

\bibitem{RonShen97} Ron A., Shen Z.,  
\emph{Affine systems in $L^2(\mathbb{R}^d)$: The analysis of the analysis operator}, 
J. Funt. Anal. 148, 408-447, 1997.

\bibitem{Ro} Rosenblatt M., 
\emph{Stationary Sequences and Random Fields}, 
Birkh\"{a}user, 1985.

\bibitem{Roz} Rozanov Y.,  
\emph{Stationary Random Processes}, 
Holden-Day, 1967.

\bibitem{samko} Samko S. G., 
\emph{Hypersingular Integrals and Their Applications}, 
CRC Press, 2002.


\bibitem{Stein} Stein E. M., 
\emph{Singular Integrals and Differentiability Properties of Functions}, 
Priceton University Press, 1987.






\end{thebibliography}
\end{document}